\documentclass{amsart}
\usepackage{amsfonts,amssymb,amsmath,amsthm}
\usepackage{url}
\usepackage{enumerate}
\usepackage{bbm}
\usepackage{mathrsfs}
\usepackage{float}
\usepackage{xypic}

\newtheorem{theo}{Theorem}[section]
\newtheorem{theoi}{Theorem}[subsection]
\newtheorem{fac}[theo]{Fact}
\newtheorem{lem}[theo]{Lemma}
\newtheorem{coro}[theo]{Corollary}
\newtheorem{propo}[theo]{Proposition}
\newtheorem{propi}[theoi]{Proposition}
\newtheorem{lemi}[theoi]{Lemma}
\newtheorem{faci}[theoi]{Fact}
\theoremstyle{definition}
\newtheorem{ttt}[theo]{}
\newtheorem{defi}[theo]{Definition}
\newtheorem{rem}[theo]{Remark}
\newtheorem{rems}[theo]{Remarks}
\newtheorem{exs}[theo]{Examples}
\newtheorem{algo}[theo]{Algorithm}
\newtheorem{ttti}[theoi]{}
\newtheorem{remi}[theoi]{Remark}
\newtheorem{remsi}[theoi]{Remarks}

\numberwithin{equation}{section}

\makeatletter
\def\hsmash{\relax 
  \ifmmode\def\next{\mathpalette\mathhsm@sh}\else\let\next\makehsm@sh
  \fi\next}
\def\makehsm@sh#1{\setbox\z@\hbox{#1}\finhsm@sh}
\def\mathhsm@sh#1#2{\setbox\z@\hbox{$\m@th#1{#2}$}\finhsm@sh}
\def\finhsm@sh{\wd\z@\z@ \box\z@}
\makeatother

\newcounter{abc}
\newenvironment{abc}{\begin{list}{\rm \alph{abc}) }{\usecounter{abc} \leftmargin=0.0pt \labelsep=0.0pt \listparindent=0.0pt \labelwidth=0.0pt \parsep=\smallskipamount \itemsep=0.0pt \topsep=0.0pt \partopsep=\smallskipamount}}{\end{list}}
\newcounter{iii}
\newenvironment{iii}{\begin{list}{\rm \roman{iii}) }{\usecounter{iii} \leftmargin=0.0pt \labelsep=0.0pt \listparindent=0.0pt \labelwidth=0.0pt \parsep=\smallskipamount \itemsep=0.0pt \topsep=0.0pt \partopsep=\smallskipamount}}{\end{list}}

\newcommand{\bbQ}{{\mathbbm Q}}
\newcommand{\bbR}{{\mathbbm R}}
\newcommand{\bbZ}{{\mathbbm Z}}

\newcommand{\Gal}{\mathop{\text{\rm Gal}}}
\newcommand{\rk}{\mathop{\text{\rm rk}}}
\newcommand{\Pic}{\mathop{\text{\rm Pic}}}
\newcommand{\Spec}{\mathop{\text{\rm Spec}}}
\newcommand{\bP}{\mathop{\text{\bf P}}\nolimits}

\newcommand{\calC}{\mathscr{C}}
\newcommand{\calI}{\mathscr{I}}
\newcommand{\calL}{\mathscr{L}}
\newcommand{\calO}{\mathscr{O}}
\newcommand{\calS}{\mathscr{S}\!}

\newcommand{\et}{{\mathop{\text{\rm {\'e}t}}}}
\newcommand{\pro}{{\rm pr}}

\newcommand{\br}{ }
\newcommand{\brr}{, }

\def\rightend#1#2{{%
 \leavevmode\nobreak\hskip .5em plus 1fil
 \penalty600 \hskip 0pt plus -1filll
 \vadjust{}\nobreak\hskip 0pt plus 1filll%
 #1\parfillskip=#2\relax \par}}

\def\eop{\ifmmode\rule[-22pt]{0pt}{1pt}\ifinner\tag*{$\square$}\else\eqno{\square}\fi\else\rightend{$\square$}{0pt}\fi}

\newcommand{\ratarrow}{$%
$\definemorphism{rat}\dashed\tip\notip%
\spreaddiagramcolumns{-12pt}%
\!\!\!\diagram%
\rrat & 
\enddiagram\!\!$%
$}

\author{Andreas-Stephan Elsenhans and J\"org Jahnel}
\address{
Universit\"at Bayreuth\\
Mathematisches Institut\\
Universit\"atsstra\ss e 30\\
D-95447 Bayreuth\\
Germany}
\email{stephan.elsenhans@uni-bayreuth.de}
\address{
Universit\"at Siegen\\
Fachbereich 6 Mathematik\\
Walter-Flex-Stra\ss e 3\\
D-57068 Siegen\\
Germany}
\email{jahnel@mathematik.uni-siegen.de}
\thanks{The computer part of this work was executed on the Sun Fire V20z Servers of the Gau\ss\ Laboratory for Scientific Computing at the G\"ottingen Mathematical Institute. Both authors are grateful to Prof.~Y.~Tschinkel for the permission to use these machines as well as to the system administrators for their~support.}

\begin{document}

\title{The discriminant of a cubic surface}

\begin{abstract}
We construct explicit examples of cubic surfaces
over~$\bbQ$
such that the 27~lines are acted upon by the index two subgroup of the maximal possible Galois~group. This~is the simple group of
order~$25\,920$.
Our~examples are given in pentahedral normal form with rational~coefficients. For~such cubic surfaces, we study the discriminant and show its relation to the index two~subgroup. On~the corresponding parameter space, we search for rational points, discuss their asymptotic, and construct an accumulating~subvariety.
\end{abstract}

\maketitle

\section{Introduction}

\begin{ttt}
Let~$\calS \subset \bP^3$
be a smooth cubic surface over an algebraically closed~field. It~is well known that there are exactly 27~lines
on~$\calS$.
The~intersection matrix of these lines is essentially the same for every smooth cubic~surface.
The~group of all permutations of the 27~lines which respect the intersection matrix is isomorphic to the Weyl
group~$W\!(E_6)$.

For a smooth cubic surface
$S \subset \bP^3$
over~$\bbQ$,
the 27~lines are, in general, not defined
over~$\bbQ$
but over an algebraic
field~extension~$L$.
The~Galois group
$\Gal (L/\bbQ)$
is a subgroup
of~$W\!(E_6)$.
It~is known that equality holds for general cubic surfaces while for diagonal cubic surfaces the Galois group is significantly~smaller. It~may be of order
$54$
at~most.
\end{ttt}

\begin{ttt}
In~this article, we describe our search for explicit examples of cubic surfaces
over~$\bbQ$
such that the~Galois group
$\Gal (L/\bbQ)$
is exactly the index~two subgroup
$D^1 W\!(E_6) \subset W\!(E_6)$.
This~is the simple group of
order~$25\,920$.

Our~approach is as~follows. We~consider cubic surfaces in pentahedral normal form with rational~coefficients. For~these, we study the
discriminant~$\Delta$.
We~show that
$\Gal (L/\bbQ)$~is
contained in the index two subgroup if and only if
$(-3)\Delta$~is
a perfect~square. This~leads to a point~search on the double covering
of~$\bP^4$
ramified at the
degree~$32$
discriminantal~variety.

A~generalized Cremona transform reduces the degree to~eight. We~discuss the asymptotic of the
\mbox{$\bbQ$-rational}
points of bounded height on the resulting double~covering and construct an accumulating~subvariety. A~final section is devoted to the problem to which extent this subvariety is~unique.
\end{ttt}

\section{The discriminant and the index two subgroup}

\begin{ttt}
One~way to write down a cubic surface in explicit form is the so-called {\em pentahedral normal form}.
Denote~by~$S^{(a_0, \ldots, a_4)}$
the cubic surface given
in~$\bP^4$
by the system of~equations
\begin{eqnarray*}
a_0 X_0^3 + a_1 X_1^3 + a_2 X_2^3 + a_3 X_3^3 + a_4 X_4^3 & = & 0 \,, \\
\phantom{a_0}\hsmash{X_0}\phantom{X_0^3} + \phantom{a_1}\hsmash{X_1}\phantom{X_1^3} + \phantom{a_2}\hsmash{X_2}\phantom{X_2^3} + \phantom{a_3}\hsmash{X_3}\phantom{X_3^3} + \phantom{a_4}\hsmash{X_4}\phantom{X_4^3} & = & 0 \,.
\end{eqnarray*}
\end{ttt}

\begin{rems}
\begin{abc}
\item
A~general cubic surface over an algebraically closed field may be brought into pentahedral normal form over that~field. Further,~the coefficients are unique up to permutation and~scaling. This~is a classical result which was first observed by J.~J.~Sylvester~\cite{Sy}. A~proof is given in~\cite{Cl}.

Cubic~surfaces in pentahedral normal form with rational coefficients are, however, special to a certain~extent.
\item
One~should keep in mind that
$S^{(0, a_1, \ldots, a_4)}$
is simply the diagonal cubic surface with coefficients
$a_1, \ldots, a_4$.
\end{abc}
\end{rems}

\begin{defi}
\label{fac}
The~expression
\begin{eqnarray*}
 & & \Delta (S^{(a_0, \ldots, a_4)}) := \\[1mm]
 & & a_0^8 \cdot\ldots\cdot a_4^8 \cdot \!\! \\
 & & \prod_{i_1, i_2, i_3, i_4 \in [0,1]} \!\!\Big( \frac1{\sqrt{a_0}} + (-1)^{i_1} \frac1{\sqrt{a_1}} + (-1)^{i_2} \frac1{\sqrt{a_2}} + (-1)^{i_3} \frac1{\sqrt{a_3}} + (-1)^{i_4} \frac1{\sqrt{a_4}} \Big)
\end{eqnarray*}
is called the {\em discriminant\/} of the cubic
surface~$S^{(a_0, \ldots, a_4)}$.
Instead~of
$\Delta (S^{(a_0, \ldots, a_4)})$,
we will usually
write~$\Delta (a_0, \ldots, a_4)$.
\end{defi}

\begin{rem}
\label{disc}
One~has
\begin{eqnarray*}
 & & \Delta (a_0, \ldots, a_4) := \\
 & & \smash{\prod_{i_1, i_2, i_3, i_4 \in [0,1]} \!\!\!\!\!\!\!\!( \sqrt{a_1 a_2 a_3 a_4} + (-1)^{i_1} \sqrt{a_0 a_2 a_3 a_4} + (-1)^{i_2} \sqrt{a_0 a_1 a_3 a_4} + \ldots} \\
 & & \hspace{5.0cm} \ldots + (-1)^{i_3} \sqrt{a_0 a_1 a_2 a_4} + (-1)^{i_4} \sqrt{a_0 a_1 a_2 a_3}) \, .
\end{eqnarray*}
\end{rem}\pagebreak[3]

\begin{lem}
$\Delta \in \bbQ[a_0, \ldots, a_4]$
is a symmetric polynomial, homogeneous of
degree~$32$,
and absolutely~irreducible.\smallskip

\noindent
{\bf Proof.}
{\em
The~remark shows
$\Delta \in \bbQ[\sqrt{a_0}, \ldots, \sqrt{a_4}]$.
Further,~the expression is obviously invariant under the action
of~$G := \Gal(\bbQ(\!\sqrt{a_0}, \ldots, \sqrt{a_4}) / \bbQ(a_0, \ldots, a_4))$.
This yields
$\Delta \in \bbQ[a_0, \ldots, a_4]$.
Symmetry~and homogeneity are~obvious.

Definition~\ref{fac} provides us with the decomposition
of~$\Delta$
into irreducible factors in the unique factorization
domain~$\smash{\overline\bbQ [\sqrt{a_0}, \ldots, \sqrt{a_4}, \frac1{a_0}, \ldots, \frac1{a_4}]}$.
Since~$G$
operates transitively on the sixteen factors, we see that
$\Delta$~is
irreducible
in~$\smash{\overline\bbQ [a_0, \ldots, a_4, \frac1{a_0}, \ldots, \frac1{a_4}]}$.

It~remains to exclude the possibility that
$\Delta$~might
be divisible by a polynomial which is a unit
in~$\smash{\overline\bbQ [a_0, \ldots, a_4, \frac1{a_0}, \ldots, \frac1{a_4}]}$.
I.e.,~by a non-trivial~monomial. For~this, note that the formula given in Remark~\ref{disc} immediately shows
$\Delta (0, a_1, a_2, a_3, a_4) = (a_1 a_2 a_3 a_4)^8$.
$\Delta$~is
not divisible
by~$a_0$.
}
\eop
\end{lem}

\begin{lem}
\label{Salmon}
Writing\/~$\sigma_i$,
for the elementary symmetric function of
degree\/~$i$
in\/~$a_0, \ldots, a_4$,
one may express the discriminant as~follows,
$$\Delta = (A^2 - 64B)^2 - 2^{11}(8D + AC) \, .$$
Here,
$$A := \sigma_4^2 - 4\sigma_3\sigma_5, \quad B := \sigma_1\sigma_5^3, \quad C := \sigma_4\sigma_5^4, \quad D := \sigma_2\sigma_5^6 \, .$$
{\bf Proof.}
{\em
This~formula may easily be established, for example, using~{\tt maple}.
}
\eop
\end{lem}

\begin{rems}
\begin{iii}
\item
Together~with
$E := \sigma_5^8$,
the expressions
$A, B, C$,
and~$D$
are called the {\em fundamental invariants\/} of the cubic
surface~$S^{(a_0, \ldots, a_4)}$.
This~notion is due to A.~Clebsch~\cite{Cl}.
\item
Lemma~\ref{Salmon}~is originally due to G.~Salmon~\cite{Sa}.
Note~that there is a misprint in Salmon's original work which has been repeatedly copied by several people throughout the 20th~century. The~correct formula may be found in~\cite{Ed}.
\end{iii}
\end{rems}

\begin{fac}
\label{sing}
Assume~that\/~$a_0 \cdot\ldots\cdot a_4 \neq 0$.
Then,~the singular points
on\/~$S^{(a_0, \ldots, a_4)}$
are exactly those of the form
$$\Big( \frac1{\sqrt{a_0}} : (-1)^{i_1} \frac1{\sqrt{a_1}} : (-1)^{i_2} \frac1{\sqrt{a_2}} : (-1)^{i_3} \frac1{\sqrt{a_3}} : (-1)^{i_4} \frac1{\sqrt{a_4}} \Big)$$
which lie on the hyperplane given
by\/~$X_0 + X_1 + X_2 + X_3 + X_4 = 0$.\smallskip

\noindent
{\bf Proof.}
{\em
A~$\overline\bbQ$-valued~point
$(x_0 : \ldots : x_4)$
on~$S^{(a_0, \ldots, a_4)}$
is singular if and only if the Jacobian~matrix
$$\left(
\begin{array}{ccccc}
3a_0x_0^2 & 3a_1x_1^2 & 3a_2x_2^2 & 3a_3x_3^2 & 3a_4x_4^2 \\
1 & 1 & 1 & 1 & 1
\end{array}
\right)$$
is not of maximal~rank. This~yields the form of the point~claimed. Observe~that
$\smash{(\frac1{\sqrt{a_0}} : (-1)^{i_1} \frac1{\sqrt{a_1}} : (-1)^{i_2} \frac1{\sqrt{a_2}} : (-1)^{i_3} \frac1{\sqrt{a_3}} : (-1)^{i_4} \frac1{\sqrt{a_4}}) \in S(\overline\bbQ)}$
if~and only if the sum of the coordinates is~zero.
}
\eop
\end{fac}

\begin{exs}
\begin{iii}
\item
The~cubic
surface~$S^{(1,1,1,1,\frac14)}$
has exactly four singular~points. These~are
$(1 : -1 : -1 : -1 : 2)$
and permutations of the first four~coordinates. This~is the famous Cayley~cubic.
\item
The~cubic
surface~$S^{(1,1,1,\frac19,\frac1{16})}$
has exactly three singular~points, namely
$(1 : -1 : -1 : -3 : 4)$
and permutations of the first three~coordinates.
\item
The~cubic
surface~$S^{(1,1,\frac14,\frac19,\frac1{25})}$
has exactly two singular~points. These~are
$(1 : -1 : -2 : -3 : 5)$
and permutations of the first two~coordinates.
\item
$(-1 : -1 : -1 : -1 : 4)$~is
the only singular point of the cubic
surface~$S^{(1,1,1,1,\frac1{16})}$.
\end{iii}
\end{exs}

\begin{coro}
\label{nonsing}
The~cubic surface\/
$S^{(a_0, \ldots, a_4)}$
is non-singular if and only if
$\Delta (a_0, \ldots, a_4) \neq 0$.\smallskip

\noindent
{\bf Proof.}
{\em
We~have that
$\Delta (0, a_1, a_2, a_3, a_4) = (a_1 a_2 a_3 a_4)^8$.
Correspondingly,~a diagonal cubic surface is singular if and only if one of its four coefficients~vanishes. In~the case
that~$a_0 \cdot\ldots\cdot a_4 \neq 0$,
the assertion follows from Fact~\ref{sing}.
}
\eop
\end{coro}

\begin{rem}
The~same is true over any ground field of
characteristic~$\neq 3$.
Therefore,~with the possible exception of the
prime~$3$,
for~$a_0, \ldots, a_4 \in \bbZ$
such
that~$\gcd (a_0, \ldots, a_4) = 1$,
the prime divisors
of~$\Delta (a_0, \ldots, a_4)$
are exactly the primes where
$S^{(a_0, \ldots, a_4)}$
has bad~reduction.

One~might want to renormalize
$\Delta$
in order to overcome the defect at the
prime~$3$.
For~this, observe that
$S^{(\frac13, \frac13, \frac13, \frac13, \frac13)}$
has only integral coefficients after the
substitution~$x_4 := -x_0 - \ldots - x_3$.
It~turns out that this surface has good reduction
at~$3$.
Since~$\Delta (\frac13, \frac13, \frac13, \frac13, \frac13) = -5 \cdot 3^{-27}$,
actually
$\pm3^{27} \Delta (a_0, \ldots, a_4)$
could have the property~desired. Theorem~\ref{klein} below indicates that the minus sign should be~correct.
\end{rem}

\begin{theo}
\label{klein}
Let\/~$a_0, \ldots, a_4 \in \bbQ$
such that\/
$\Delta (a_0, \ldots, a_4) \neq 0$.
Then,~the Galois group operating on the 27~lines
on\/~$S^{(a_0, \ldots, a_4)}$
is contained in the index two
subgroup\/~$D^1 W\!(E_6) \subset W\!(E_6)$
if and only if\/
$(-3) \Delta (a_0, \ldots, a_4) \in \bbQ$
is a perfect~square.\smallskip

\noindent
{\bf Proof.}
{\em
{\em First step.}
Construction of a ramified covering of degree two
of~$\bP^4$.\smallskip

\noindent
Define~$\smash{\calC \subset \bP^4_{(X)} \times \bP^4_{(x)}}$
by the system of equations
\begin{eqnarray*}
x_0 X_0^3 + x_1 X_1^3 + x_2 X_2^3 + x_3 X_3^3 + x_4 X_4^3 & = & 0 \,, \\
\phantom{x_0}\hsmash{X_0}\phantom{X_0^3} + \phantom{x_1}\hsmash{X_1}\phantom{X_1^3} + \phantom{x_2}\hsmash{X_2}\phantom{X_2^3} + \phantom{x_3}\hsmash{X_3}\phantom{X_3^3} + \phantom{x_4}\hsmash{X_4}\phantom{X_4^3} & = & 0 \,.
\end{eqnarray*}
The~projection
$\smash{\pi \colon \calC \to \bP^4 (= \bP^4_{(x)}})$
is the family of the cubic surfaces in pentahedral normal~form. The~fiber
of~$\pi$
over~$(x_0 : \ldots : x_4)$
is the cubic
surface~$S^{(x_0, \ldots , x_4)}$.

The~fiber~$\calC_\eta$
over the generic
point~$\eta \in \bP^4$
is a smooth cubic surface
over~$\bbQ(\eta) = \bbQ(x_1/x_0, x_2/x_0, x_3/x_0, x_4/x_0)$.
Its~27~lines are defined over a finite
extension~$L$
of~$\bbQ(\eta)$.
We~claim
that~$\Gal (L/\bbQ(\eta)) = W\!(E_6)$.

Indeed,~this is the maximal possible~group. The~inclusion
``$\subseteq$''
is, therefore, trivially~fulfilled. On~the other hand, according to a result of B.~L.~van~der~Waerden, the generic Galois group
$\Gal (L/\bbQ(\eta))$
can not be smaller than that for a particular~fiber. Specializing,~for example, to
$(x_0 : x_1 : x_2 : x_3 : x_4) = (1 : 2 : 3 : 7 : 17)$,
\mbox{\cite[Algorithm~10]{EJ}} shows that the Galois group is equal
to~$W\!(E_6)$.

Consequently,~there exists a unique intermediate
field~$K$
of~$L/\bbQ(\eta)$
which is quadratic
over~$\bbQ(\eta)$.
This~induces a
scheme~$V$
together with a finite morphism
$p\colon V \to \bP^4$
of degree~two.

In~fact, this is a standard~construction. For~each affine open set
$\Spec A = U \subseteq \bP^4$,
take the spectrum of the integral closure
of~$A$
in the
extension~$K$.
Note~that
$A$~is
integrally closed
in~$\bbQ(\eta)$
since
$\bP^4$
is a normal~scheme. The~morphism
$p\colon V \to \bP^4$
is finite according to the finiteness of the integral~closure.\medskip\pagebreak[3]

\noindent
{\em Second step.}
$p \colon V \to \bP^4$~is
unramified outside the
divisor~$R$
given
by~``$\Delta = 0$''.\smallskip

\noindent
For~this, let us describe the double
covering~$V$
more~precisely. We have
$\bP^3_{\bbQ(\eta)} \subset \bP^4_{\bbQ(\eta)}$
given by the equation
$X_0 + \ldots + X_4 = 0$
and a smooth cubic
surface
$$\calC_\eta \subset \bP^3_{\bbQ(\eta)} \, .$$
On~$\calC_\eta$,
there are the 45~tritangent~planes. These~give rise to a subscheme
of~$\smash{(\bP^3)^\vee_{\bbQ(\eta)}}$
which is finite of length~45 and \'etale
over~$\bbQ(\eta)$.

This,~according to Galois theory, induces a
set~$M = \{e_1, \ldots , e_{45}\}$
of~45~elements together with an operation
of~$\Gal(\overline{Q(\eta)}/Q(\eta))$.
Actually,~only a finite quotient isomorphic
to~$W\!(E_6)$
is~operating.
The~set~$M$,
in turn, gives rise to the two element
set~$\{\pm e_1 \wedge ... \wedge e_{45}\}$
which is again acted upon
by~$\Gal(\overline{Q(\eta)}/Q(\eta))$.
The~fixgroup of this operation corresponds to the quadratic field
extension~$K/\bbQ(\eta)$.

The~same may be done in the relative situation
over~$\bP^4 \setminus R$.
The~45 tritangent~planes yield a closed subscheme
of~$(\bP^3)^\vee \times \bP^4$
which is finite and \'etale of degree~45
over~$\bP^4 \setminus R$.
According~to A.~Grothendieck's theory of the \'etale fundamental group~\cite{SGA1}, this induces a
set~$M = \{e_1, \ldots , e_{45}\}$
of~45~elements together with an operation of~$\pi_1^\et (\bP^4 \setminus R, *)$.
This~group is canonically a quotient
of~$\Gal(\overline{Q(\eta)}/Q(\eta))$.
Again,~we get a canonical operation on the two element
set~$\{\pm e_1 \wedge ... \wedge e_{45}\}$.
Corresponding~to this, there is an \'etale
covering~$p^\prime\colon V^\prime \to \bP^4 \setminus R$
of degree~two.

$V^\prime$~is,
by construction, a normal scheme with function
field~$K$.
In~particular, over an affine open set
$\Spec A = U \subseteq \bP^4 \setminus R$,
we have the spectrum of the integral closure
of~$A$
in the
extension~$K$.
This~shows that
$V$
and~$V^\prime$
coincide
over~$\bP^4 \setminus R$.\medskip

\noindent
{\em Third step.}
The~equation.\smallskip

\noindent
As~$R$
is irreducible, the ramification locus 
of~$p\colon V \to \bP^4$
might be either empty or equal
to~$R$.
If~the ramification locus were empty then,
as~$\pi_1^\et (\bP^4, *) = 0$,
we had a trivial covering by a non-connected~scheme.
However,~$V$
is connected by~construction. The~generic fiber
of~$p$
is a scheme consisting of a single~point.

Hence,~$p$
is ramified exactly
at~$R$.
This~implies
that~$V$
is given by the~equation
$w^2 = \lambda \Delta$
for a suitable
constant~$\lambda$.\medskip\pagebreak[3]

\noindent
{\em Fourth step.}
Specialization.\smallskip

\noindent
Let~$(a_0 : \ldots : a_4) \in \bP^4 (\bbQ)$
such
that~$\Delta (a_0, \ldots, a_4) \neq 0$.
Then,~by virtue of the construction above, we have the following~statement.

Denote~by~$l$
the field of definition of the 27~lines
on~$S^{(a_0, \ldots, a_4)}$.
Then,~the smallest intermediate
field~$k$
of~$l/\bbQ$
such that
$\Gal (l/k)$~acts
on the 45~tritangent~planes
on~$S^{(a_0, \ldots, a_4)}$
only via even permutations is exactly
$k = \bbQ (\sqrt{\lambda \Delta (a_0, \ldots, a_4)})$.

This~extension splits if and only
if~$\lambda \Delta (a_0, \ldots, a_4)$
is a perfect square
in~$\bbQ$.
Except~for the determination of the
constant~$\lambda$,
this proves the~assertion.\medskip

\noindent
{\em Fifth step.}
The~constant~$\lambda$.\smallskip

\noindent
We~consider the particular cubic
surface~$S^{(0,1,1,1,1)}$.
I.e.,~the diagonal cubic surface given
by~$x_1^3 + x_2^3 + x_3^3 + x_4^3 = 0$.

Here,~the 27~lines are defined over the
field~$\bbQ(\zeta_3) = \bbQ(\sqrt{-3})$.
They~may be given explicitly in the form
$$x_i + \zeta_3^m x_j = 0, \qquad x_k + \zeta_3^n x_l = 0$$
for~$\{ i, j, k, l \} = \{ 1, 2, 3, 4 \}$
and~$m,n \in \{0,1,2\}$.
Exactly~three of these lines are defined
over~$\bbQ$.

They~form a triangle which is cut out by the
equation~$x_1 + x_2 + x_3 + x_4 = 0$.
There~are six more tritangent~planes which consist of a rational line and two lines conjugate to each~other. These~are given
by~$x_i + x_j = 0$
for~$\{i,j\} \subset \{1, 2, 3, 4\}$
any subset of size~two. To~summarize,
$\Gal(\bbQ(\sqrt{-3})/\bbQ)$
operates on the 45~tritangent~planes as a product of 19~two-cycles while seven tritangent~planes are~fixed. This~is an odd~permutation.

Consequently,~in this case,
$k = \bbQ (\sqrt{-3})$
is the smallest field such that
$\Gal (l/k)$~acts
on the 45~tritangent~planes only by even~permutations.
As~$\Delta (0,1,1,1,1) = 1$,
this
shows~$\lambda = -3$
up to a factor which is a perfect~square. The~proof is~complete.
}
\eop
\end{theo}

\begin{rem}
This~result was essentially known to H.~Burk\-hardt~\cite[p.~341]{Bu} in~1893. Burkhardt~gives credit to C.~Jordan~\cite{Jo} who was the first to study the automorphism group of the configuration of the 27~lines on a cubic~surface.
\end{rem}

\section{Rational points on the discriminantal covering}

\begin{defi}
We~will call the twofold covering
of~$\bP^4_\bbQ$
given by the~equation
\begin{equation}
\label{Zwei}
w^2 = -3\Delta(a_0, \ldots, a_4)
\end{equation}
the {\em discriminantal~covering}.
\end{defi}

\begin{ttt}
There are two surprising constraints which equation~(\ref{Zwei}) imposes on the coefficients~$a_0, \ldots, a_4$.
\end{ttt}\pagebreak[3]

\begin{propo}[{\rm The two constraints}{}]
~\!\!\!--\!--\!--\!--\;
Suppose\/~$a_0, \ldots, a_4 \in \bbZ$
are such
that\/~$\gcd (a_0, \ldots, a_4) = 1$
and\/~$(-3) \Delta(a_0, \ldots, a_4) \neq 0$
is a perfect square
in\/~$\bbQ$.

\begin{abc}
\item
Then,~$a_0, \ldots, a_4$
all have the same~sign.
\item
Further,~for every prime
number\/
$p \equiv 2 \pmod 3$,
all the\/
$p$-adic~valuations
$\nu_p (a_0), \ldots, \nu_p (a_4)$
are~even.
\end{abc}\smallskip

\noindent
{\bf Proof.}
{\em
Observe~first that the assumption
ensures~$a_0, \ldots, a_4 \neq 0$.
Indeed,
$\Delta (0, a_1, \ldots, a_4) = (a_1a_2a_3a_4)^8 \geq 0$.\smallskip

\noindent
a)
Assume~the~contrary. Then,~there is a product of four, say
$a_1 \cdot\ldots\cdot a_4$,
which is~negative. The~formula given in Remark~\ref{disc} implies that
$\Delta(a_0, \ldots, a_4)$
is the norm of an element
of~$\bbQ (\sqrt{a_1 \cdot\ldots\cdot a_4})$.
As~this is an imaginary quadratic~field, we see
that~$\Delta(a_0, \ldots, a_4) \geq 0$.
Contradiction!\smallskip

\noindent
b)
Again,~assume the~contrary. Then,~there is a product of four, say
$a_1 \cdot\ldots\cdot a_4$,
the
$p$-adic~valuation
of which is~odd. We~have the fact that
$\Delta(a_0, \ldots, a_4)$
is the norm of an element
of~$\bbQ (\sqrt{a_1 \cdot\ldots\cdot a_4})$.
On~the other hand,
$(-3)\Delta(a_0, \ldots, a_4)$,
being perfect square by assumption, is a norm,~too.
Consequently,~$(-3)$
is the norm of an element
of~$\bbQ (\sqrt{a_1 \cdot\ldots\cdot a_4})$.

Since~$\nu_p (a_1 \cdot\ldots\cdot a_4)$
is odd, the norm~equation
$(-3) = x^2 - a_1 \cdot\ldots\cdot a_4 \cdot y^2$
ensures that
$\nu_p (x) = 0$
and~$\nu_p (a_1 \cdot\ldots\cdot a_4 \cdot y^2) > 0$.
Therefore,~$(-3)$
is a quadratic residue
modulo~$p$.
This~is a~contradiction.
}
\eop
\end{propo}

\begin{ttt}
We~are interested in smooth cubic
surfaces~$S^{(a_0, \ldots, a_4)}$
such that the Galois group operating on the 27~lines is exactly equal
to~$D^1 W\!(E_6)$.

By~Theorem~\ref{klein}, this implies that
$(a_0 : \ldots : a_4) \in \bP^4 (\bbQ)$
gives rise to a
\mbox{$\bbQ$-rational}~point
on the discriminantal~covering. Further,~according to Corollary~\ref{nonsing},
$(a_0 : \ldots : a_4)$~is
supposed not to lie on the ramification~locus.

Finally,~if two of the coefficients were the same,
say~$a_0 = a_1$,
then
$S^{(a_0, \ldots, a_4)}$
allowed the tritangent plane
``$x_0 + x_1 = 0$''
which was defined
over~$\bbQ$.
Consequently,~the order of the group acting on the lines could not be higher
than~$1152$.
\end{ttt}

\begin{ttt}
{\bf A naive search.}
For~these reasons, we searched for
\mbox{$\bbQ$-r}a\-tio\-nal points
$(w; \,a_0 : \ldots : a_4)$
satisfying equation~(\ref{Zwei}) and the extra conditions~below,

\begin{iii}
\item
$w \neq 0$.
\item
No~two of the five
coordinates~$a_0, \ldots, a_4$
are the~same.
\end{iii}

\noindent
A~rather simple computation led to the
\mbox{$\bbQ$-rational}
points
$(3 : 4 : 21 : 36 : 63)$,
$(4 : 7 : 12 : 28 : 84)$,
and~$(12 : 28 : 36 : 63 : 84)$.
Up~to~symmetry, these are the only solutions of
height~$\leq 100$.
\end{ttt}

\begin{rem}
The~three rational points given above really lead to cubic surfaces such that the 27~lines are acted upon by the simple
group~$D^1 W\!(E_6)$.
To~prove this, we ran the algorithm below which is an obvious modification of~\mbox{\cite[Algorithm~10]{EJ}}.
\end{rem}

\begin{algo}[{\rm Verifying
$G \supseteq D^1 W\!(E_6)$}{}]
\label{WE6}
~\!--\!--\!--\!--\;
Given the
equation~$f=0$
of a smooth cubic surface, this algorithm verifies that
$G \subseteq W\!(E_6)$
is of index at most~two.

\begin{iii}
\item
Compute~a univariate polynomial
$0 \neq g \in \bbZ[d]$
of minimal degree such that
$$g \in (f(\ell(0)), f(\ell(\infty)), f(\ell(1)), f(\ell(-1))) \subset \bbQ[a,b,c,d]$$
where~$\ell \colon t \mapsto (1 : t : (a + b t) : (c + d t)).$

\noindent
If~$g$
is not of
degree~$27$
then terminate with an error~message. In~this case, the coordinate system is not sufficiently~general.
\item
Factor~$g$
modulo all primes below a given~limit. Ignore~the primes dividing the leading coefficient
of~$g$.
\item\label{deco}
If~one of the factors is multiple then go to the next prime~immediately.
Otherwise, check whether the decomposition type
is~$(1,1,5,5,5,5,5)$
or~$(9,9,9)$.
\item
If~each of the two cases occurred at least once then output the message \mbox{``The~Galois group contains~$D^1 W\!(E_6)$.''}~and~terminate.

\noindent
Otherwise,~output ``Can not prove that the Galois group contains~$D^1 W\!(E_6)$.''
\end{iii}
\end{algo}

\section{The generalized Cremona transform}

\begin{ttt}
$\Delta$~is
a homogeneous form of
degree~$32$.
Naively,~one would expect that there are not many solutions of the~equation
$$w^2 = -3\Delta(a_0, a_1, a_2, a_3, a_4) \, .$$
The~constraints proven above reduce expectations even~more. Nevertheless,~three rational points of
height~$\leq \!100$
have been~found. The~reason for this is the following~observation.
\end{ttt}\pagebreak[3]

\begin{fac}
There~is
form\/~$\Delta^\prime$
homogeneous of
degree\/~$8$
such~that
$$\Delta (a_0, \ldots, a_4) = (a_0 \cdot\ldots\cdot a_4)^8 \cdot \Delta^\prime (1/a_0, \ldots, 1/a_4) \, .$$
{\bf Proof.}
{\em
The~octic~$\Delta^\prime$
is given by the~formula
\begin{eqnarray*}
 & & \Delta^\prime (x_0, \ldots, x_4) := \\
 & & \hspace{1cm} \smash{\!\!\!\!\!\!\prod_{i_1, i_2, i_3, i_4 \in [0,1]} \!\!\!\!\!\!\!\!\!\!\big( \sqrt{x_0} + (-1)^{i_1} \sqrt{x_1} + (-1)^{i_2} \sqrt{x_2} + (-1)^{i_3} \sqrt{x_3} + (-1)^{i_4} \sqrt{x_4} \big) \, .}
\end{eqnarray*}
}\vskip-3mm
\eop
\end{fac}

\begin{defi}
We~will call the birational
automorphism~$\iota$
of~$\bP^4$
given~by
$$(a_0 : \ldots : a_4) \mapsto (1/a_0 : \ldots : 1/a_4)$$
a {\em generalized Cremona transform}. Note~that the standard Cremona transform
of~$\bP^2$
is given
by~$(a_0 : a_1 : a_2) \mapsto (1/a_0 : 1/a_1 : 1/a_2)$.\smallskip

\noindent
The~generalized Cremona
transform~$\iota$
provides
a bijection~of
$$\{ (x_0 : \ldots : x_4) \in \bP^4 (\bbQ) \mid x_0 \cdot\ldots\cdot x_4 \neq 0 \}$$
to~itself.
\end{defi}

\begin{coro}
$(x_0 : \ldots : x_4) \in \bP^4 (\bbQ)$,
$x_0 \cdot\ldots\cdot x_4 \neq 0$,
gives rise to a solution~of
$$w^2 = (-3) \Delta^\prime (x_0, \ldots, x_4)$$
if and only if\/
$\iota \big( (x_0 : \ldots : x_4) \big)$
yields a rational point on the discriminantal~covering.
\eop
\end{coro}

\begin{lem}
\label{quadrat}
\begin{abc}
\item
$\Delta^\prime \in \bbQ[x_0, \ldots, x_4]$
is a symmetric polynomial, homogeneous of
degree eight and absolutely~irreducible.
\item
One~has\/
$\smash{\Delta^\prime (0, x_1, \ldots, x_4) = D^2}$
for a symmetric, homogeneous quartic form\/
$D \in \bbQ[x_1, \ldots, x_4]$.
\end{abc}\medskip

\noindent
{\bf Proof.}
{\em
a)
By~definition,
$\Delta^\prime \in \bbQ[\sqrt{x_0}, \ldots, \sqrt{x_4}]$.
Further,~the expression
for~$\Delta^\prime$
is obviously invariant under the action
of~$G := \Gal(\bbQ(\!\sqrt{x_0}, \ldots, \sqrt{x_4}) / \bbQ(x_0, \ldots, x_4))$.
This~yields
$\Delta \in \bbQ[x_0, \ldots, x_4]$.
Symmetry~and homogeneity are~obvious.

Finally,~we have a decomposition
of~$\Delta^\prime$
into irreducible factors in the unique factorization
domain~$\smash{\overline\bbQ [\sqrt{x_0}, \ldots, \sqrt{x_4}]}$.
Since~$G$~operates
transitively on the sixteen factors,
$\Delta$~is
absolutely~irreducible.\smallskip

\noindent
b)
$\Delta^\prime (0, x_1, \ldots, x_4)$~is
the square~of
$$D (x_1, \ldots, x_4) := \!\!\!\!\!\!\prod\limits_{i_2, i_3, i_4 \in [0,1]} \!\!\!\!\!\!\!\!\big( \!\sqrt{x_1} + (-1)^{i_2} \sqrt{x_2} + (-1)^{i_3} \sqrt{x_3} + (-1)^{i_4} \sqrt{x_4} \big) \, . \eop$$
}
\end{lem}

\begin{rems}
\begin{iii}
\item
The~ramification locus
$R := \text{''}\!\Delta^\prime\! = 0$''
is a rational~threefold. The~parametrization
$\iota\colon \bP^3 \to R$
given~by
$$\iota\colon (t_0 : \ldots : t_3) \mapsto \big( t_0^2 : t_1^2 : t_2^2 : t_3^2 : (t_0 + \ldots + t_3)^2 \big)$$
is a finite birational~morphism.
\item
The~equation
$D = 0$
defines the Roman surface of~J.~Steiner.
\end{iii}
\end{rems}

\section{More rational points on the discriminantal covering}

\subsection{A point search}

\begin{ttti}
On~the double
covering~$\pi\colon O \to \bP^4_\bbQ$,
given~by
$$w^2 = (-3) \Delta^\prime (x_0, \ldots, x_4) \, ,$$
we searched for rational points such~that

\begin{iii}
\item
$w \neq 0$,
\item
the five
coordinates~$x_0, \ldots, x_4$
are pairwise different from each~other.
\end{iii}
%
\end{ttti}

\begin{ttti}
Surprisingly~many solutions have been~found. It~turned out that there are
$4\,900\,907$
essentially different solutions up to a height limit
of~$3000$.
Under~symmetry, they give rise to
$120$~solutions~each.
The~smallest ones are
$(1 : 3 : 7 : 9 : 12)$,
$(1 : 3 : 4 : 7 : 13)$,
$(1 : 3 : 7 : 12 : 13)$,
and~$(3 : 7 : 9 : 12 : 13)$.
For~a few height limits, we indicate the number of solutions up to that limit in the table~below.

\begin{table}[H]
\caption{Numbers of solutions up to various height limits}
\centerline{
\begin{tabular}{|r|r||r|r||r|r||r|r|}
\hline
 limit & \# & limit & \# & limit & \# & limit & \# \\
\hline
  25 &     20 &  200 & 10\,039 &   500 &  93\,680 &  1500 & 1\,111\,303 \\
  50 &    209 &  300 & 25\,778 &   750 & 236\,403 &  2000 & 2\,088\,752 \\
 100 & 1\,481 &  400 & 54\,331 &  1000 & 460\,330 &  3000 & 4\,900\,907 \\
\hline
\end{tabular}}
\label{Tab1}
\end{table}
\end{ttti}

\begin{remi}
We~used the constraints shown above to optimize the searching~algorithm. On~one hand, it is sufficient to search for solutions such that
$0 < x_0 < x_1 < x_2 < x_3 < x_4$.
On~the other hand, only
$751$
of the positive integers up
to~$3000$
fulfill the condition that all prime divisors
$p \equiv 2 \pmod 3$
have an even~exponent.
\end{remi}

\subsection{The conjecture of Manin}

\begin{ttti}
Let~$X$
be a non-singular (weak) Fano~variety
over~$\bbQ$.
Assume~that
$X(\bbQ) \neq \emptyset$.
Then,~the conjecture of Manin~\cite{FMT} makes the following prediction for the number of
\mbox{$\bbQ$-rational}
points
on~$X$
of bounded anticanonical~height.

There~exists
some~$\tau > 0$
such that, for every Zariski open
set~$X^\circ \subseteq X$
which is sufficiently small but non-empty,
$$\# \{ x \in X^\circ (\bbQ) \mid h_{-K} (x) < B \} \,\sim\, \tau B \log^r \!B$$
for~$r := \rk \Pic (X) - 1$
and~$B \gg 0$.

Unfortunately,~$O$
is~singular. In~this situation, one has to consider a
resolution~$\widetilde{O}$
of singularities and compare~heights.
\end{ttti}

\begin{propi}
The~singular locus
of\/~$O$
is reducible into ten components. The~component\/
$S_{(x_0, x_1)}$
is given~by
$$x_0 - x_1 = 0 \, , \quad x_2^2 + x_3^2 + x_4^2 - 2x_2x_3 - 2x_2x_4 - 2x_3x_4 = 0 \, .$$
The~others are obtained by permuting~coordinates.\medskip

\noindent
{\bf Proof.}
{\em
{\em First case.}
$x_0 \cdot\ldots\cdot x_4 \neq 0$.\smallskip

\noindent
Then,~the morphism
$p \colon \bP^4_\bbQ \to \bP^4_\bbQ$
given~by
$$(t_0 : \ldots : t_4) \mapsto (t_0^2 : \ldots : t_4^2)$$
is \'etale
over~$(x_0 : \ldots : x_4)$.
We~may therefore test the fiber~product
$O \times_{\pi, \bP^4_\bbQ, p} \bP^4_\bbQ$
for~smoothness. It~is given explicitly~by
$$w^2 = (-3) \!\!\!\!\!\!\prod_{i_1, i_2, i_3, i_4 \in [0,1]} \!\!\!\!\!\!\!\!\big( t_0 + (-1)^{i_1} t_1 + (-1)^{i_2} t_2 + (-1)^{i_3} t_3 + (-1)^{i_4} t_4 \big) \, .$$
Here,~the singular points are exactly the singular points of the ramification~locus. That,~in turn, consists of 16~hyperplanes such that precisely the intersection points are~singular. Going~back
to~$O$,
we see that the singular points are those where at least two of the~expressions
$$\sqrt{x_0} + (-1)^{i_1} \sqrt{x_1} + (-1)^{i_2} \sqrt{x_2} + (-1)^{i_3} \sqrt{x_3} + (-1)^{i_4} \sqrt{x_4}$$
vanish.

If~these expressions coincide in one or four signs then this enforces one coordinate to be~zero. The~cases that there are two or three signs in common are essentially equivalent to each~other. Without~restriction,
$$\sqrt{x_0} - \sqrt{x_1} + \sqrt{x_2} + \sqrt{x_3} + \sqrt{x_4} = \sqrt{x_0} - \sqrt{x_1} - \sqrt{x_2} - \sqrt{x_3} - \sqrt{x_4} = 0 \, .$$
Then,~$\sqrt{x_0} = \sqrt{x_1}$
and~$\sqrt{x_2} + \sqrt{x_3} + \sqrt{x_4} = 0$.
The~first equation
yields~$x_0 = x_1$.
The~quadratic relation given is equivalent to
$\sqrt{x_2} \pm \sqrt{x_3} \pm \sqrt{x_4} = 0$.\medskip\pagebreak[3]

\noindent
{\em Second case.}
$x_0 \cdot\ldots\cdot x_4 = 0$.\smallskip

\noindent
The~singular locus is a Zariski closed~subset. Therefore,~the points satisfying the equations given above are clearly~singular. It~remains to prove that the others are~non-singular.

Without~restriction, we may assume that
$x_0 = 0$
and that exactly one of the
expressions~$\sqrt{x_1} + (-1)^{i_2} \sqrt{x_2} + (-1)^{i_3} \sqrt{x_3} + (-1)^{i_4} \sqrt{x_4}$,
say
$\sqrt{x_1} + \sqrt{x_2} + \sqrt{x_3} + \sqrt{x_4}$,
is equal to~zero. Then,~the partial derivative~of
\begin{eqnarray*}
(\sqrt{x_0} + \sqrt{x_1} + \sqrt{x_2} + \sqrt{x_3} + \sqrt{x_4})(\sqrt{x_0} - \sqrt{x_1} - \sqrt{x_2} - \sqrt{x_3} - \sqrt{x_4}) = {}\hspace{0.45cm} \\
{} = x_0 - (\sqrt{x_1} + \sqrt{x_2} + \sqrt{x_3} + \sqrt{x_4})^2
\end{eqnarray*}
by~$x_0$
is non-zero. As~the other factors do not vanish, the product over all the 16~factors has non-zero derivative at this~point. The~assertion~follows.
}
\eop
\end{propi}

\begin{theoi}
\label{resol}
Let\/~$\pro\colon \widetilde{O} \to O$
be the proper and birational morphism obtained by blowing up the ten singular~components.

\begin{abc}
\item
Then,~$\widetilde{O}$
is non-singular.
I.e.,~$\pro$
is a resolution of~singularities.
\item
Further,~$\rk \Pic (\widetilde{O}) = 11$.
\item
The~canonical divisor
of\/~$\smash{\widetilde{O}}$
is\/~$K = \pro^* K_O$
for\/
$K_O = - \pi^* H$
and\/
$H$
a hyperplane~section
of\/~$\bP^4$.
\end{abc}\smallskip

\noindent
{\bf Proof.}
{\em
a)
This~may be tested~locally.
Let~$(w; \,x_0 : \ldots : x_4)$
be a point in the singular locus
of~$O$.\medskip\pagebreak[3]

\noindent
{\em First case.}
$x_0 \cdot\ldots\cdot x_4 \neq 0$.\smallskip

\noindent
Near~$(x_0 : \ldots : x_4)$,
the~morphism
$$p \colon \bP^4_\bbQ \to \bP^4_\bbQ, \quad (t_0 : \ldots : t_4) \mapsto (t_0^2 : \ldots : t_4^2)$$
is~\'etale. \;We~may take square roots
$\smash{t_0^{(0)}, \ldots, t_4^{(0)}}$
of~$x_0, \ldots, x_4$
and~consider
$$\smash{w^2 = (-3) \!\!\!\!\!\!\!\prod_{i_1, i_2, i_3, i_4 \in [0,1]} \!\!\!\!\!\!\!\!\!\big( t_0 + (-1)^{i_1} t_1 + (-1)^{i_2} t_2 + (-1)^{i_3} t_3 + (-1)^{i_4} t_4 \big) \, .}$$
Actually,~only the linear factors vanishing
at~$\smash{(t_0^{(0)} : \ldots : t_4^{(0)})}$
need to be taken into~consideration.

Without~restriction, suppose
that~$(x_0 : \ldots : x_4) \in S_{(x_0, x_1)}$.
Then,~again without~restriction,
$$\smash{t_0^{(0)} - t_1^{(0)} + t_2^{(0)} + t_3^{(0)} + t_4^{(0)} = t_0^{(0)} - t_1^{(0)} - t_2^{(0)} - t_3^{(0)} - t_4^{(0)}} = 0 \, .$$
The~corresponding linear
forms~$X, Y$
are linearly independent which means that we blow up a scheme, locally given by the equation
$W^2 = XY$,
at the
ideal~$(X, Y)$.
The~result is clearly non-singular.\smallskip

Now~suppose that
$(x_0 : \ldots : x_4)$
is a point of intersection of at least two singular components. Without~loss of generality, the second singular component might be either
$S_{(x_0, x_2)}$
or~$S_{(x_2, x_3)}$.
The~latter variant enforces that
$(x_0 : \ldots : x_4) = (1 : 1 : 1 : 1 : 4)$
is the point corresponding to the Cayley~cubic. This~is actually a special case of the first~variant.

Thus,~assume that\vspace{0.4mm}
$(x_0 : \ldots : x_4) \in S_{(x_0, x_1)} \cap S_{(x_0, x_2)}$.
Then,~without \mbox{restriction,}
$\smash{t_0^{(0)} = t_1^{(0)} = t_2^{(0)}}$
and~$\smash{t_0^{(0)} + t_3^{(0)} + t_4^{(0)} = 0}$.
We~have the three vanishing linear forms
$t_0 + t_1 - t_2 + t_3 + t_4$,
$t_0 - t_1 + t_2 + t_3 + t_4$,
and~$t_0 - t_1 - t_2 - t_3 - t_4$.
Only~when
$x_3 = x_0$
(or~$x_4 = x_0$),
another linear form~vanishes.

Altogether,~there are four linearly independent linear
forms~$X, Y, Z$,
and~$U$.
We~blow
up~$W^2 = XY\!ZU$
or~$W^2 = XY\!Z$
at~$(X, Y)$,
$(X, Z)$,
and~$(Y, Z)$,
(as well
as~$(X, U)$,
$(Y, U)$,
and~$(Z, U)$).
The~resulting scheme is non-singular.\medskip\pagebreak[3]

\noindent
{\em Second case.}
Exactly one of the coordinates
$x_0, \ldots, x_4$~vanishes.\smallskip

\noindent
Then,~without loss of generality,
$(x_0 : \ldots : x_4) = (a : a : b : b : 0)$.
We~may take square roots
$t_0, \ldots, t_3$
of~$x_0, \ldots, x_3$
such that
$t_0$
and~$t_1$
as well as
$t_2$
and~$t_3$
are of the same~sign. Then,~the right hand side goes over into the product over
all~$(t_0 \pm t_1 \pm t_2 \pm t_3)^2 - x_4$.
Among~these,
$(t_0 - t_1 + t_2 - t_3)^2 - x_4$
and~$(t_0 - t_1 - t_2 + t_3)^2 - x_4$
do~vanish. 

Hence,~for two linearly independent linear
forms
$X$
and~$Y$,
we consider the scheme given
by~$W^2 = (X^2 - x_4)(Y^2 - x_4)$.
The~singular components
$S_{(x_0, x_1)}$
and~$S_{(x_2, x_3)}$
correspond to the ideals
$(X^2 - x_4, X + Y)$
and~$(X^2 - x_4, X - Y)$,~respectively.
Blowing~up the first ideal amounts to the substitutions
$x_4 = X^2 + v(X + Y)$
and, for the other affine chart,
$x_4 = X^2 + \frac1v (X + Y)$.
The~first substitution leads to
$(W^\prime)^2 = v(X - Y + v)$
becoming smooth after blowing~up
$(v, X - Y)$
which is the next~step. On~the other hand, the second substitution yields
$(W^\prime)^2 = v(X - Y) + 1$
which is clearly non-singular
near~$v = 0$.\smallskip

There~is the exceptional case
that~$a = b$.
Then,~$(t_0 + t_1 - t_2 - t_3)^2 - x_4$
is a third factor~vanishing. We~have to consider a scheme locally given by
$W^2 = (X^2 - x_4)(Y^2 - x_4)(Z^2 - x_4)$.
Here,~the substitution
$x_4 = X^2 + v(X + Y)$
yields~$(W^\prime)^2 = v(X - Y + v)[Z^2 - X^2 - v(X + Y)]$.
The~next step, to blow up
$(v, X - Y)$,
leads to
$(W^{\prime\prime})^2 = v_1(1 + v_1)[Z^2 - X^2 - v_1(X^2 - Y^2)]$.
Here,~for the other affine chart, we find a formula of the same~structure. Further,~it is sufficient to consider the singularity
at~$v_1 = 0$.
That~at~$v_1 = -1$
is~analogous.

Actually,~to blow up\vspace{0.2mm}
$S_{(x_1, x_2)} \cup S_{(x_0, x_3)}$
suffices to resolve this~singularity. Indeed,~the substitution
$Z^2 - X^2 = v_2v_1$
yields
$(W^{\prime\prime\prime})^2 = v_2 - X^2 + Y^2$
which is clearly non-singular. On~the other hand, putting
$\smash{Z^2 - X^2 = \frac1{v_2} v_1}$
leads to
$(W^{\prime\prime\prime})^2 = v_2 (1 - v_2(X^2 - Y^2))$
which is obviously smooth
near~$v_2 = 0$.\medskip\pagebreak[3]

\noindent
{\em Third case.}
Exactly two of the coordinates
$x_0, \ldots, x_4$~vanish.\smallskip

\noindent
Here,~without restriction,
$(x_0 : \ldots : x_4) = \big( 0 : 0 : \big(t_2^{(0)}\big)^2 : \big(t_3^{(0)}\big)^2 : \big(t_4^{(0)}\big)^2 \big)$
for~$\smash{t_2^{(0)} + t_3^{(0)} + t_4^{(0)} = 0}$.
Therefore,~precisely four of the sixteen factors of the right hand side~vanish. These~are
$\sqrt{x_0} \pm \sqrt{x_1} \pm (t_2 + t_3 + t_4)$.
We~find
$W^2 = X^2 - 2YT^2 + T^4$
for the new coordinate functions
$X := x_0 - x_1$,
$Y := x_0 + x_1$,
and~$T := t_2 + t_3 + t_4$.

When~blowing up
$(X, T)$,
the substitution
$X := uT$
leads to
$(W^\prime)^2 = u^2 - 2Y + T^2$
which is non-singular. On~the other hand,
$T := uX$
yields
$(W^\prime)^2 = 1 - 2u^2Y + u^4X^2$
being clearly smooth
near~$u = 0$.\medskip\pagebreak[3]

\noindent
{\em Fourth case.}
Three~of the coordinates
$x_0, \ldots, x_4$~vanish.\smallskip

\noindent
Without~restriction,
$(x_0 : \ldots : x_4) = (0 : 0 : 0 : 1 : 1)$.
Take~square roots
$t_3, t_4$
of~$x_3$
and~$x_4$
which are of the same~sign. The~eight factors
$\sqrt{x_0} \pm \sqrt{x_1} \pm \sqrt{x_2} \pm (t_3 - t_4)$
vanish
at~$(0 : 0 : 0 : 1 : 1)$.
We~find the local equation
$W^2 = D(x_0, x_1, x_2, t^2)$
for
$t := t_3 - t_4$
and
$D$~the
symmetric, homogeneous quartic from Lemma~\ref{quadrat}.b).

Blowing~up
$(x_0 - x_1, x_2 - t^2)$
amounts to substituting
$x_2 := t^2 + u(x_0 - x_1)$
and, for the other affine chart,
$x_2 := t^2 + \frac1u (x_0 - x_1)$.
Then,~the ideals
$(x_0 - x_2, x_1 - t^2)$
and
$(x_1 - x_2, x_0 - t^2)$
to be blown up subsequently go over to
$(u - 1, x_1 - t^2)$
and~$(u + 1, x_0 - t^2)$.
The~substitutions
\begin{eqnarray*}
x_2 & := & t^2 + u(x_0 - x_1) \, , \\
x_1 & := & t^2 + u_1(u - 1) \, , \\
x_0 & := & t^2 + u_2(u + 1) \, ,
\end{eqnarray*}
yield
$$W^2 = (Y - uZ)^2 + 8Zt^2$$
for the new functions
$Y := u_1 + u_2$
and~$Z := u_1 - u_2$.

For~the other seven affine charts of this triple blow-up, the equations are completely~analogous. The~differences are that the definitions of
$Y$
and~$Z$
may be replaced
by~$Y,Z := 1 \pm u_1u_2$.
Further,~instead of
$Y - uZ$,
we may
have~$uY - Z$.

The~last step is to blow up the
ideal~$(Y - uZ, t)$
corresponding to the
component~$S_{(x_3, x_4)}$.
The~substitution
$Y - uZ = vt$
yields
$(W^\prime)^2 = v^2 + 8Z$
which is non-singular. Indeed,~otherwise we must have
$v = 0$
and~$W^\prime = 0$
which
implies~$Z = 0$.
But,~in this situation,
$Z$~is
a local~parameter. On~the other hand,
$Y - uZ = \frac1v t$
leads to
$(W^\prime)^2 = 1 + 8v^2 Z$
which is clearly smooth
at~$v = 0$.\smallskip

\noindent
b) We~claim that
$O$~is~normal.
To~see this, note first
that~$O$
is a hypersurface in weighted projective
space~$\bP = \bP (4,1,1,1,1,1)$.
This~is a scheme equipped with a canonical rational
map~$\iota\colon \bP \;\ratarrow\, \bP (1,1,1,1,1) = \bP^4$.
$\iota$~is
undefined at exactly one point which is the only singularity
of~$\bP$.

By~construction,~the double
covering~$O$
does not meet the singular~point.
Consequently,~$O$
is Gorenstein and, in particular, Cohen-Macaulay. Further,~the singularities
of~$O$
are in
codimension~$2$.
Serre's~criterion~\cite[Theorem 23.8]{Ma} shows that
$O$~is~normal.

We~assert that, after each step of blowing~up, the resulting scheme is still~normal. In~fact, the centre of the blowing~up is a codimension~two complete~intersection.
The~blow-up~$\smash{{\rm Bl}_{S_{(x_0, x_1)}} (O)}$
is, therefore, locally given by a single equation in
a~\mbox{$\bP^1$-bundle}
over~$O$.
This~ensures
${\rm Bl}_{S_{(x_0, x_1)}} (O)$~is
Cohen-Macaulay. Further,~the smooth part
of~$O$
is untouched under~blowing~up. Thus,~regularity in codimension two could be destroyed only if the whole exceptional
set were~singular. As~this is a
\mbox{$\bP^1$-bundle}
over~$S_{(x_0, x_1)}$,
that is clearly not the~case. The~same argument works for each of the subsequent~steps.

By~Lemma~\ref{inj}, it suffices to show that the Picard rank grows by one in each~step. Again,~let us explain this for the first step in order to simplify~notation. We~have
${\rm Bl}_{S_{(x_0, x_1)}} (O) = {\bf Proj} (\calO \oplus \calI \oplus \calI^2 \oplus \ldots \,)$
for~$\calI := \calI_{S_{(x_0, x_1)}, O}$.
We~assert that the twisting
sheaf~$\calO(1)$
is linearly independent of the pull-backs
of~$\Pic (O)$
in~$\smash{\Pic ({\rm Bl}_{S_{(x_0, x_1)}} (O))}$.
Indeed,~$\calO (n)$
for~$n \neq 0$
is non-trivial when restricted to one of the exceptional fibers which is just
a~$\bP^1$.\smallskip

\noindent
c)
As~$O$
is a Gorenstein scheme, its dualizing
sheaf~$\omega_O$
is invertible~\cite[Theorem~3.5.1]{Co}. To~describe
$\omega_O$
completely, we may restrict it
to~$O^{\rm reg}$
since
$O$
is~normal. Here,~$\omega_O |_{O^{\rm reg}} \cong \Omega_{O^{\rm reg}}^4$.
A~$4$-form
with a simple pole
at~``$x_0 = 0$''
is given
by~$({x_0^4}/w) \cdot d(x_1/x_0) \wedge \ldots \wedge d(x_4/x_0)$.
Hence,~$\omega_O = \pi^* \calO (-1)$.

Further,~$\pro$
is an isomorphism outside the exceptional~fibers. This~implies that
$K$
and~$\pro^* K_O$
coincide up to a sum of exceptional~divisors. Due~to symmetry, the coefficients at~$E_1, \ldots, E_{10}$
are equal to each~other. To~determine the actual number, consider a general point~$P \in S_{(x_0, x_1)}$.
Near~$P$,
we blow up a double covering of the
type~$w^2 = XY$.
This~is a quadric cone times a neighbourhood
of~$(0, 0) \in {\bf A}^2$.
Its~blow-up is the Hirzebruch
surface~$\Sigma_2$
times that~neighbourhood. The~exceptional curve
$E \subset \Sigma_2$
is
a~\mbox{$(-2)$-curve},
hence
$\omega_{\Sigma_2} |_E$
is~trivial. The~coefficients desired are equal to~zero.
}
\eop
\end{theoi}

\begin{lemi}
\label{inj}
Let\/~$p\colon X \to Y$
be a surjective and birational morphism of Noetherian, normal, integral~schemes. Then,~the pull-back homomorphism
$p^* \colon \Pic(Y) \to \Pic(X)$
is~injective.\smallskip

\noindent
{\bf Proof.}
{\em
Suppose,~for
$\calL \in \Pic(Y)$,
the pull-back
$p^* \calL \in \Pic(X)$
would be~trivial. This~means, we have a section
$s \in \Gamma (X, p^* \calL)$
without zeroes or~poles. Corresponding~to each codimension one
point~$\xi \in Y$,
there is a discrete valuation
ring~$\calO_\xi$.
Further,~there is a codimension one
point~$\zeta \in X$
mapping
to~$\xi$.
As~$\smash{\calO_\xi}$
is integrally closed, we see
that~$\smash{\calO_\xi \cong \calO_\zeta}$.

Consequently,~$s$~gives rise
to a section
$t \in \Gamma (Y^\circ, \calL |_{Y^\circ})$
without zeroes or poles
for~$Y^\circ \subseteq Y$
the complement of a closed subset of
codimension~$\geq 2$.
\cite[Theorem 12.4.i)]{Ma}~implies that
$t$~may
be extended to a global~section.
Hence,~$\calL \cong \calO_Y$
is~trivial.
}
\eop
\end{lemi}

\begin{remi}[{\rm The prediction---Manin's conjecture for the double covering $O$}{}]
~\\
Theorem~\ref{resol}.c)~implies
$$h_{-K} (y) = h_{-\pro^* K_O} (y) = h_{-K_O} (\pro(y)) = h_{{\rm naive}, \bP^4} (\pi ( \pro(y)))$$
for
every~$\smash{y \in \widetilde{O} (\bbQ)}$.
Manin's~conjecture therefore predicts that, for every sufficiently small, non-empty, Zariski open
subset~$O^\circ \subseteq O$,
$$\# \{ x \in O^\circ (\bbQ) \mid h_{{\rm naive}, \bP^4} (\pi(x)) < B \} \,\sim\, \tau B \log^{10} \!B \, .$$
The~reader might want to compare Table~\ref{Tab2}~below where the actual numbers are given for a reasonably chosen Zariski open~subset.
%
\end{remi}

\begin{remi}
We~actually found that
$\Pic (\widetilde{O}) \cong \bbZ^{11}$
is a trivial
$\Gal(\overline\bbQ/\bbQ)$-module.
This~implies that there is no Brauer-Manin obstruction present
on~$\smash{\widetilde{O}}$.
\end{remi}\pagebreak[3]

\subsection{Infinitely many solutions}
\label{Infi}

\begin{propi}
\label{inf}
There~are infinitely many\/
$\bbQ$-rational
points
on\/~$O$.
In~fact, over the quadric
surface\/~$Q$
in\/~$\bP^4$,
given~by\/~$l = q = 0$~for\pagebreak[3]
\begin{eqnarray*}
l & := & x_0 + x_1 + x_2 - 3x_3 - 3x_4 \, , \\
q & := & x_0^2 + x_1^2 + x_2^2 + 9x_3^2 - x_0x_1 - x_0x_2 - 3x_0x_3 - x_1x_2 - 3x_1x_3 - 3x_2x_3 \, ,
\end{eqnarray*}
the double covering\/
$\pi \colon O \to \bP^4_\bbQ$~splits.
In~particular, there are one or two\/
\mbox{$\bbQ$-rational}
points above each\/
\mbox{$\bbQ$-rational}
point
of\/~$Q$.\smallskip

\noindent
{\bf Proof.}
{\em
Modulo~$\calI_Q$,
one has~actually
\begin{equation}
\label{Drei}
\textstyle{(-3) \Delta(x_0, \ldots, x_4) = [\frac{64}3 (x_0 - x_1)(x_0 - x_2)(x_1 - x_2)(x_3 - x_4)]^2 \, .}
\end{equation}
}\vskip-6.5mm
\eop\vskip4mm
\end{propi}

\begin{remsi}
\begin{iii}
\item
The~difference of the two octic forms in equation~(\ref{Drei}) consists
of~$495$~monomials.
To~verify the assertion, one may first use the linear equation to
eliminate~$x_4$
and then check that the remaining octic form
in~$x_0, \ldots, x_3$
is divisible by the quadratic
form~$q$.

Actually,~a simple Gr\"obner base calculation quarries the fact that equation~(\ref{Drei}) is true even
modulo~$\calI_Q^2$.
\item
There~is another~proof for Lemma~\ref{inf} which is somehow easier from the computational point of view but less canonical.
In~fact,
$Q$~is
parametrized
by~the birational~map
$\iota \colon \bP^2 \,\ratarrow\, Q$,
\begin{eqnarray*}
            (t_0 : t_1 : t_2) \mapsto \hspace{9.8cm} \\
        \big( (t_0^2 + t_1^2 + t_2^2 - t_0t_1  -t_0t_2 - t_1t_2)
            : (t_0^2 + t_1^2 + t_2^2 - t_0t_1 +2t_0t_2 - t_1t_2) :
                                                            \hspace{0.1cm} \\
            : (t_0^2 + t_1^2 + t_2^2 - t_0t_1  -t_0t_2 +2t_1t_2)
            : t_2^2
            : (t_0^2 + t_1^2 - t_0t_1) \big)
\end{eqnarray*}
which is defined
over~$\bbQ$.
The~locus where
$\iota$
is undefined does not contain any
\mbox{$\bbQ$-rational}
point since the quadratic form
$t_0^2 + t_1^2 - t_0t_1$
does not represent zero
over~$\bbQ$.
A~direct calculation~shows
$$(-3) \Delta^\prime \big( \iota(t_0, t_1, t_2) \big) = [576 t_0 t_1 (t_0 - t_1) t_2^3 (t_0^2 + t_1^2 - t_0t_1 - t_2^2)]^2 \, .$$
Here,~the
factor~$t_0$
corresponds
to~$(x_0 - x_1)$,
$t_1$
to~$(x_0 - x_2)$,
$(t_0 - t_1)$
to~$(x_1 - x_2)$,
and~$(t_0^2 + t_1^2 - t_0t_1 - t_2^2)$
to~$(x_3 - x_4)$.
The~factor~$t_2^3$
is somehow~artificial.
For~$t_2 = 0$,
the parametrization is constant
to~$(1 : 1 : 1 : 0 : 1)$.

The~parametrization~$\iota$
is actually constructed in a very naive~manner. Start~with the point
$(1 : 1 : 1 : 0 : 1)$
and determine for which value
of~$\tau \neq 0$
the~point
$$(1 : (1 + \tau t_0) : (1 + \tau t_1) : (\tau t_2/3) : (1 + \tau (t_0 + t_1 - t_2)/3))$$
is contained in the quadric
surface~$Q$.
Many~other parametrizations would serve the same~purpose.
\end{iii}
\end{remsi}

\begin{remsi}
\begin{iii}
\item
The~surface~$Q$
is obviously symmetric under permutations
of~$\{x_0, x_1, x_2\}$.
It~is symmetric under switch
of~$x_3$
and~$x_4$,~too.
All~in all, there are ten mutually different copies 
of~$Q$.
\item
$Q$~is
a smooth quadric~surface. The~two pencils of lines
on~$Q$
are defined
over~$\bbQ(\sqrt{-3})$
and conjugate to each~other.
\item
This~implies
that~$\Pic (Q) = \bbZ$.
The~Picard group has two generators as soon as the ground field
contains~$\bbQ(\sqrt{-3})$.

For~quadrics such
as~$Q$,
Manin's conjecture is~proven. The~number of points of
height~$\leq B$
is
asymptotically~$\tau_Q B^2$
for some
constant~$\tau_Q$.
This~means that
$\pi^{-1} (Q) \subset O$
is an example of a so-called {\em accumulating~subvariety.} The~growth of the number of rational points
on~$\pi^{-1} (Q)$
is faster than predicted for a sufficiently small Zariski open subset
of~$O$.
\end{iii}
\end{remsi}

\begin{remi}
The~quadric surface
$Q$~was
detected by a statistical investigation of the rational points found
on~$O$.
Nevertheless,~as the height limit
of~$3000$
is too low, most of these points are actually not contained
in~$\pi^{-1} (Q)$
or one of its~copies. Cf.~Table~\ref{Tab2} below for the numbers of points
on~$O$
with those over the copies
of~$Q$~excluded.

\begin{table}[H]
\caption{Numbers of solutions, accumulating subvarieties excluded}
\centerline{
\begin{tabular}{|r|r||r|r||r|r||r|r|}
\hline
 limit & \# & limit & \# & limit & \# & limit & \# \\
\hline
  25 &     12 &  200 &  8\,989 &   500 &  86\,897 &  1500 & 1\,049\,502 \\
  50 &    156 &  300 & 23\,496 &   750 & 221\,187 &  2000 & 1\,977\,863 \\
 100 & 1\,248 &  400 & 50\,070 &  1000 & 432\,737 &  3000 & 4\,651\,857 \\
\hline
\end{tabular}}
\label{Tab2}
\end{table}
\end{remi}

\begin{remsi}
\begin{iii}
\item
The~smallest
$\bbQ$-rational points
on~$Q$
with no two coordinates equal
are~$(3 : 9 : 12 : 1 : 7)$
and~$(1 : 7 : 13 : 3 : 4)$.
Algorithm~\ref{WE6} shows that, indeed, these two points yield cubic surfaces such that the 27~lines are acted upon by the simple
group~$D^1 W\!(E_6)$.
According~to B.~L.~van~der~Waerden, this is the generic behaviour
on~$Q$.
\item
When~testing the cubic surface corresponding to
$(3 : 9 : 12 : 1 : 7)$,
Algorithm~\ref{WE6} works with the primes
$19$
and~$73$.
Therefore,~we have an explicit infinite set of
\mbox{$\bbQ$-rational} points
which lead to the
group~$D^1 W\!(E_6)$.
It~is given by those points
on~$Q$
reducing
to~$(3 : 9 : 12 : 1 : 7)$
modulo
both~$19$
and~$73$.
\end{iii}
\end{remsi}

\begin{ttti}
Some~of the surprising properties
of~$Q$
are described by the following two~facts.
\end{ttti}

\begin{faci}
$Q$~meets the
octic~$R$
only within its singular~locus. Actually,
$$Q \cap R \,\subset\, S_{(x_0, x_1)} \cup S_{(x_0, x_2)} \cup S_{(x_1, x_2)} \cup S_{(x_3, x_4)} \, .$$%
\noindent
{\bf Proof.}
{\em
Suppose~$(x_0 : \ldots : x_4) \in Q \cap R$.
Then,~formula~(\ref{Drei})
implies that
$x_0 = x_1$,
$x_0 = x_2$,
$x_1 = x_1$,
or~$x_3 = x_4$.
The~equation~$x_0 = x_1$
yields
$x_0 = (-x_2 + 3x_3 + 3x_4)/2$.
Substituting~this into the quadratic
relation~$q(x_0, \ldots, x_4) = 0$
from the definition
of~$Q$~shows
$$x_2^2 + x_3^2 + x_4^2 - 2x_2x_3 - 2x_2x_4 - 2x_3x_4 = 0 \, .$$
For~the~relations
$x_0 = x_2$,
$x_1 = x_2$,
and~$x_3 = x_4$,
the situation is~analogous.
}
\eop
\end{faci}

\begin{faci}
$Q$~is
tangent to all five coordinate~hyperplanes.\smallskip

\noindent
The~points of tangency
are\/~$(0 : 3 : 3 : 1 : 1)$,
$(3 : 0 : 3 : 1 : 1)$,
$(3 : 3 : 0 : 1 : 1)$,
$(1 : 1 : 1 : 0 : 1)$,
and~$(1 : 1 : 1 : 1 : 0)$.
\eop
\end{faci}

\begin{ttti}
The~quadric
surface~$Q$
determines the linear
form~$l$~uniquely.
On~the other hand, the quadratic
form~$q$
is unique only up to a multiple
of~$l$.
One~might have the idea to fix a canonical
representative~$\underline{q}$
by the requirement that the quadric
threefold~``$\underline{q} = 0$''
contain some of the singular components~entirely. This~is possible to a certain~extent.\medskip

\noindent
{\bf Fact.}
{\em
{\rm a)}
There~is no quadric threefold
in\/~$\bP^4$
containing the singular components\/
$S_{(x_0, x_1)}$
and\/~$S_{(x_3, x_4)}$.\smallskip

\noindent
{\rm b)}
There~is, however, a one-dimensional family of quadric threefolds
in\/~$\bP^4$
containing\/
$S_{(x_0, x_1)}$
and\/~$S_{(x_0, x_2)}$.
It~is given by\/
$f_t = 0$
for a
parameter\/~$t$~and
\begin{eqnarray*}
& & f_t := -t x_0^2 + x_1^2 + x_2^2 + x_3^2 + x_4^2 - {} \\
& &\hspace{1.3cm} {} - (1 - t)x_0x_1 - (1 - t)x_0x_2 + 2x_0x_3 + 2x_0x_4 + {} \\
& &\hspace{1.6cm} {} + (1 - t)x_1x_2 - 2x_1x_3 - 2x_1x_4 - 2x_2x_3 - 2x_2x_4 - 2x_3x_4  = 0 \, .
\end{eqnarray*}%
}%
\noindent
{\bf Proof.}
The~statement that a quadric threefold
contains~$S_{(x_0, x_1)}$
is equivalent to saying it is given by an equation of the
form~$\underline{q} = 0$~for
$$\underline{q} := a (x_2^2 + x_3^2 + x_4^2 - 2x_2x_3 - 2x_2x_4 - 2x_3x_4) + (a_0x_0 + \ldots + a_4x_4) \cdot (x_0 - x_1) \, .$$
The~assumptions of~a) yield a linear system of equations which is only trivially~solvable. On~the other hand, the system of equations for~b) leads to a two-dimensional vector~space.
\eop
\end{ttti}

\begin{remi}
This~family is attached to the rational~map
$f \colon \bP^4 \,\ratarrow\, \bP^1$,
\begin{eqnarray*}
 & & \hspace{-0.2cm} (x_0 : \ldots : x_4) \mapsto (x_1^2 + x_2^2 + x_3^2 + x_4^2 - x_0x_1 - x_0x_2 + 2x_0x_3 + 2x_0x_4 + {} \\
 & & \hspace{3.0cm} {} + x_1x_2 - 2x_1x_3 - 2x_1x_4- 2x_2x_3 - 2x_2x_4 - 2x_3x_4) \\
 & & \hspace{7.5cm} {} : (x_0^2 - x_0x_1 - x_0x_2 + x_1x_2) \, .
\end{eqnarray*}
The~map
$f$
enjoys the following remarkable~properties.

\begin{iii}
\item
Its~locus of indeterminacy is equal
to~$S_{(x_0, x_1)} \cup S_{(x_0, x_2)}$.
\item
The~fiber
at~$t = -1$
is a singular quadric of rank~three. The~fiber at infinity is reducible into the two
hyperplanes~``$x_0 = x_1$''
and~``$x_0 = x_2$''.
All~other special fibers are~smooth.
\item
The~special fiber
at~$t = \frac13$
may also be written as
$$4q + (-7x_0 + 5x_1 + 5x_2 + 9x_3 - 3x_4)l = 0 \, .$$
In~particular, the accumulating
subvariety~$Q$
is contained within this~fiber.
\item
The~fiber
at~$t = \frac13$
contains more of the rational points known than any~other, even after deleting the accumulating~subvarieties. The~singular fiber
at~$t = -1$
follows~next.
\end{iii}
\end{remi}

\section{Accumulating subvarieties}

\begin{ttt}
The~goal of this section is to prove that there are no other accumulating subvarieties which are, in a certain sense, similar
to~$Q$.
Similarity~shall include to be a non-degenerate quadric surface over which the double
covering~$\pi\colon O \to \bP^4_\bbQ$~splits.

In~view of the first constraint established above, this implies that the real points on such a
quadric~surface~$S$
are contained in the~16-ant
$$\{\, (x_0 : \ldots : x_4) \in \bP^4(\bbR) \mid x_0, \ldots, x_4 \geq 0 {\rm ~or~} x_0, \ldots, x_4 \leq 0 \,\} \, .$$
Further,~there are strong restrictions for the behaviour at the~boundary. By~Lemma~\ref{quadrat}.b), we know that
$\Delta^\prime$~is
a perfect square on the coordinate
hyperplane~$H_0$
given
by~``$x_0 = 0$''.
On~the other hand, we require
$(-3) \Delta^\prime$
to be a perfect square
on~$S$.

A~way to realize both of these, seemingly contradictory, requirements is to
make~$S \cap H_0$
a curve of degree two on which
$(-3)$~is
the square of a rational~function. The~only such examples are two lines
over~$\bbQ(\sqrt{-3})$
which are conjugate to each~other. This~implies that
$S$
must necessarily be tangent
to~$H_0$
and the point of tangency is a
$\bbQ$-rational~point
on~the ramification
locus~$R$.
\end{ttt}

\begin{theo}
\label{main}
Suppose\/~$\smash{S \subset \bP^4_\bbQ}$
is a smooth quadric surface such that the double
covering\/~$\pi\colon O \to \bP^4_\bbQ$
splits
over\/~$S$.
Assume~further that\/
$S$~is
tangent to the five coordinate hyperplanes\/
$H_0, \ldots, H_4$
and that, for
each\/~$i$,
the point of tangency is actually contained in one of the three lines
on\/~$H_i \cap R$.\smallskip

\noindent
Then,~$S$
is equal
to\/~$Q$
or one of its copies under permutation of~coordinates.
\end{theo}

\begin{rem}
On~the Steiner
surface~$H_0 \cap R$,
there are two types of
$\bbQ$-ra\-tio\-nal~points.
There~are the three lines given
by~$(0 : r : r : s : s)$
and permutations of the four coordinates to the~right. The~other $\bbQ$-ra\-tio\-nal
points are of the form
$(0 : t_1^2 : \ldots : t_4^2)$
for~$t_1, \ldots, t_4 \in \bbQ$
such
that~$t_1 + \ldots + t_4 = 0$.
\end{rem}

\begin{lem}
Assume\/~$S$
is as in Theorem~\ref{main}. Further,~write
$$P^{(0)} := (0 : x_1^{(0)} : x_2^{(0)} : x_3^{(0)} : x_4^{(0)})$$
for the point of tangency
of\/~$S$
with the coordinate
hyperplane\/~$H_0$.\smallskip

\noindent
Then,~$\smash{x_1^{(0)}, x_2^{(0)}, x_3^{(0)}, x_4^{(0)} \neq 0}$.\medskip

\noindent
{\bf Proof.}
{\em
Assume, to the contrary,
that~$x_1^{(0)} = 0$.
The~assumption on the type of the points of tangency made in Theorem~\ref{main} implies that one more coordinate must~vanish. Without~restriction, we may
assume~$P^{(0)} = (0 : 0 : 0 : 1 : 1)$.
The~tangent plane
at~$P^{(0)}$
is given by
$x_0 = 0$
and another linear
relation~$C_1 x_1 + \ldots + C_4 x_4 = 0$.
Whatever~the coefficients are, there is a tangent vector
$(v_0, \ldots, v_4)$
such that
$v_1 < 0$
or~$v_2 < 0$.
The~implicit function theorem yields a real point
$(x_0 : \ldots : x_3 : 1) \in S(\bbR)$
satisfying
$x_1 < 0$
or~$x_2 < 0$.
This~is a~contradiction.
}
\eop
\end{lem}

\begin{lem}
\label{Groeb}
Assume that the quadric
surface\/~$S$
is tangent to the coordinate
hyperplanes\/~$H_0$,\vspace{0.5mm}
$H_1$,
and\/~$H_2$
in\/~$\smash{(0 : x_1^{(0)} : x_2^{(0)} : x_3^{(0)} : x_4^{(0)})}$,
$\smash{(x_0^{(1)} : 0 : x_2^{(1)} : x_3^{(1)} : x_4^{(1)})}$,
and\/~$\smash{(x_0^{(2)} : x_1^{(2)} : 0 : x_3^{(2)} : x_4^{(2)})}$,~respectively.\smallskip

\noindent
Then,
\begin{equation*}
\hspace{2.3cm} x_1^{(0)} x_2^{(1)} x_0^{(2)} - x_2^{(0)} x_0^{(1)} x_1^{(2)} = 0
\end{equation*}
or
\begin{eqnarray*}
x_1^{(0)} x_2^{(1)} x_0^{(2)} + x_2^{(0)} x_0^{(1)} x_1^{(2)} & = & 0 \, , \nonumber \\
x_1^{(0)} x_0^{(1)} x_3^{(2)} - x_1^{(0)} x_0^{(2)} x_3^{(1)} - x_0^{(1)} x_1^{(2)} x_3^{(0)} & = & 0 \, , \nonumber \\
x_1^{(0)} x_2^{(1)} x_3^{(2)} + x_2^{(0)} x_1^{(2)} x_3^{(1)} - x_2^{(1)} x_1^{(2)} x_3^{(0)} & = & 0 \, , \\
x_2^{(0)} x_0^{(1)} x_3^{(2)} - x_2^{(0)} x_0^{(2)} x_3^{(1)} + x_2^{(1)} x_0^{(2)} x_3^{(0)} & = & 0 \, . \nonumber
\end{eqnarray*}
{\bf Proof.}
{\em
The~linear equation by which
$S$~is
defined may be written
\begin{equation}
\label{1}
L_0x_0 + L_1x_1 + L_2x_2 + L_3x_3 + L_4x_4 = 0 \, .
\end{equation}
We~distinguish three~cases.\medskip\pagebreak[3]

\noindent
{\em First case.}
$L_4 \neq 0$.\smallskip

\noindent
Then,~we may use the linear equation~(\ref{1}) to
eliminate~$x_4$
from the quadratic~equation. Write
\begin{eqnarray*}
 Q_0x_0^2 + Q_1x_1^2 + Q_2x_2^2 + Q_3x_3^2 + {} \hspace{6.4cm} \\ 
 {} + Q_4x_0x_1 + Q_5x_0x_2 + Q_6x_0x_3 + Q_7x_1x_2 + Q_8x_1x_3 + Q_9x_2x_3 & = & 0 \, . \nonumber
\end{eqnarray*}
Tangency~of
$H_0$~at
$(0 : x_1^{(0)} : x_2^{(0)} : x_3^{(0)} : x_4^{(0)})$
means that the two linear~forms
\begin{eqnarray*}
 (Q_4x_1^{(0)}+Q_5x_2^{(0)}+Q_6x_3^{(0)}) x_0 + (2Q_1x_1^{(0)}+Q_7x_2^{(0)}+Q_8x_3^{(0)}) x_1 + {} \hspace{1.05cm} \\
\hspace{1.05cm} {} + (2Q_2x_2^{(0)}+Q_7x_1^{(0)}+Q_9x_3^{(0)}) x_2 + (2Q_3x_3^{(0)}+Q_8x_1^{(0)}+Q_9x_2^{(0)}) x_3 \, , \\
 L_0x_0 + L_1x_1 + L_2x_2 + L_3x_3 + x_4 \, , \,
\end{eqnarray*}
together
generate~$x_0$.
This~enforces the linear~relations
\begin{eqnarray}
\label{tang_rel}
  2x_1^{(0)} Q_1 + x_2^{(0)} Q_7 + x_3^{(0)} Q_8 & = & 0 \, , \nonumber \\
  2x_2^{(0)} Q_2 + x_1^{(0)} Q_7 + x_3^{(0)} Q_9 & = & 0 \, ,           \\
  2x_3^{(0)} Q_3 + x_1^{(0)} Q_8 + x_1^{(0)} Q_9 & = & 0 \, . \nonumber 
\end{eqnarray}
The~two other points of tangency yield relations which are completely~analogous. Altogether,~we find the homogeneous linear system of equations associated with
the~\mbox{$9 \times 10$-matrix}
$$\left(
\begin{array}{cccccccccc}
     0 &  2x_1^{(0)} &    0 &    0 &    0 &    0 &    0 &   x_2^{(0)} &   x_3^{(0)} &    0 \\
     0 &    0 &  2x_2^{(0)} &    0 &    0 &    0 &    0 &   x_1^{(0)} &    0 &   x_3^{(0)} \\
     0 &    0 &    0 &  2x_3^{(0)} &    0 &    0 &    0 &    0 &   x_1^{(0)} &   x_2^{(0)} \\
   2x_0^{(1)} &    0 &    0 &    0 &    0 &   x_2^{(1)} &   x_3^{(1)} &    0 &    0 &    0 \\
     0 &    0 &  2x_2^{(1)} &    0 &    0 &   x_0^{(1)} &    0 &    0 &    0 &   x_3^{(1)} \\
     0 &    0 &    0 &  2x_3^{(1)} &    0 &    0 &   x_0^{(1)} &    0 &    0 &   x_2^{(1)} \\
   2x_0^{(2)} &    0 &    0 &    0 &   x_1^{(2)} &    0 &   x_3^{(2)} &    0 &    0 &    0 \\
     0 &  2x_1^{(2)} &    0 &    0 &   x_0^{(2)} &    0 &    0 &    0 &   x_3^{(2)} &    0 \\
     0 &    0 &    0 &  2x_3^{(2)} &    0 &    0 &   x_0^{(2)} &    0 &   x_1^{(2)} &    0
\end{array}
\right)
.
$$
If~this matrix is of
rank~$9$
then the quadratic equation
defining~$S$
is, up to scaling, determined~uniquely. In~fact, this case is~degenerate. There~is a linear form
in~$x_0, \ldots, x_3$
only, vanishing on the three points~given. The~unique solution of the system corresponds to the square of this linear~form.

Consequently,~the rank is at
most~$8$.
The~ten
\mbox{$9 \times 9$-minors}
must all~vanish. These~minors are polynomials
in~$\smash{x_0^{(0)}, \ldots, x_3^{(2)}}$
having
$$(x_1^{(0)} x_2^{(1)} x_0^{(2)} - x_2^{(0)} x_0^{(1)} x_1^{(2)})$$
as their greatest common~divisor. After~division by this, we are left with~ten~\mbox{sextics}. It~turns out that they are precisely the squares and pairwise products of the~four~cubics
\mbox{$\smash{x_1^{(0)} x_2^{(1)} x_0^{(2)} + x_2^{(0)} x_0^{(1)} x_1^{(2)}}$,
$x_1^{(0)} x_0^{(1)} x_3^{(2)} - x_1^{(0)} x_0^{(2)} x_3^{(1)} - x_0^{(1)} x_1^{(2)} x_3^{(0)}$,\vspace{0.6mm}}
$x_1^{(0)} x_2^{(1)} x_3^{(2)} + x_2^{(0)} x_1^{(2)} x_3^{(1)} - x_2^{(1)} x_1^{(2)} x_3^{(0)}$,
and
$\smash{x_2^{(0)} x_0^{(1)} x_3^{(2)} \!- x_2^{(0)} x_0^{(2)} x_3^{(1)} \!+ x_2^{(1)} x_0^{(2)} x_3^{(0)}}$.\medskip

\noindent
{\em Second case.}
$L_4 = 0$
and~$L_3 \neq 0$.\smallskip

\noindent
As~the roles of the third and fourth coordinates may be interchanged, we have, as~in the first case,
$\smash{x_1^{(0)} x_2^{(1)} x_0^{(2)} - x_2^{(0)} x_0^{(1)} x_1^{(2)} = 0}$
or
$$x_1^{(0)} x_2^{(1)} x_0^{(2)} + x_2^{(0)} x_0^{(1)} x_1^{(2)} = 0 \, .$$
Suppose~that the second variant is~present. Then,~the linear equation~(\ref{1}) implies that the vector\vspace{0.5mm}
$\smash{(x_3^{(0)}, x_3^{(1)}, x_3^{(2)})^t}$
is linearly dependent of
$\smash{(0, x_0^{(1)}, x_0^{(2)})^t}$,
$\smash{(x_1^{(0)}, 0, x_1^{(2)})^t}$,
and
$\smash{(x_2^{(0)}, x_2^{(1)}, 0)^t}$.\vspace{0.2mm}
For~these vectors instead of
$\smash{(x_3^{(0)}, x_3^{(1)}, x_3^{(2)})^t}$,
the three more relations asserted are clearly~true.\medskip

\noindent
{\em Third case.}
$L_3 = L_4 = 0$.\smallskip

\noindent
In~this situation, we may write the three points of tangency in the form
\mbox{$\smash{(0 : L_2 :\! (-L_1) \!:\! x_3^{(0)} \!:\! x_4^{(0)})}$,\vspace{0.5mm}
$\smash{(L_2 : 0 :\! (-L_0) \!:\! x_3^{(1)} \!:\! x_4^{(1)})}$,
and
$\smash{(L_1 :\! (-L_0) : 0 \!:\! x_3^{(2)} \!:\! x_4^{(2)})}$.}
It~turns out that the relation
$$x_1^{(0)} x_2^{(1)} x_0^{(2)} + x_2^{(0)} x_0^{(1)} x_1^{(2)} = 0$$
is automatically~fulfilled.
Further,~$L_0, L_1, L_2 \neq 0$.
Each~of the three equations still~to be proven reduces
to~$\smash{L_0x_3^{(0)} - L_1x_3^{(1)} + L_2x_3^{(2)} = 0}$.

We~may use the linear equation~(\ref{1}) to eliminate
$x_0$
from the quad\-ratic~equation. Write
\begin{eqnarray*}
 Q_0x_1^2 + Q_1x_2^2 + Q_2x_3^2 + Q_3x_4^2 + {} \hspace{6.4cm} \\ 
 {} + Q_4x_1x_2 + Q_5x_1x_3 + Q_6x_1x_4 + Q_7x_2x_3 + Q_8x_2x_4 + Q_9x_3x_4 & = & 0 \, . \nonumber
\end{eqnarray*}
Tangency~of
$H_0$
at~$(0 : L_2 : (-L_1) : x_3^{(0)} : x_4^{(0)})$
yields the linear~relations\pagebreak[3]
\begin{eqnarray*}
L_2(2L_2Q_0 - L_1Q_4 + x_3^{(0)}Q_5 + x_4^{(0)}Q_6) - {} \hspace{4.9cm} \\
{} - L_1(-2L_1Q_1 + L_2Q_4 + x_3^{(0)}Q_7 + x_4^{(0)}Q_8) & = & 0 \, , \\
2x_3^{(0)}Q_2 + x_1^{(0)}Q_5 + x_2^{(0)}Q_7 + x_4^{(0)}Q_9 & = & 0 \, , \\
2x_4^{(0)}Q_3 + x_1^{(0)}Q_6 + x_2^{(0)}Q_8 + x_3^{(0)}Q_9 & = & 0 \, .
\end{eqnarray*}
Tangency~of
$H_1$
and~$H_2$
leads to linear relations completely analogous to those given in~(\ref{tang_rel}). Altogether,~we find the homogeneous linear system of equations associated with
the~\mbox{$9 \times 10$-matrix}
$$\left(
\begin{array}{cccccccccc}
2L_2^2 & 2L_1^2 & 0 & 0 & -2L_1\!L_2 & x_3^{(0)}\!L_2 & x_4^{(0)}\!L_2 & -x_3^{(0)}\!L_1 & -x_4^{(0)}\!L_1 & 0 \\
   0 &     0 & 2x_3^{(0)} &           0 &    0 & x_1^{(0)} &         0 & x_2^{(0)} &         0 & x_4^{(0)} \\
   0 &     0 &          0 &  2x_4^{(0)} &    0 &         0 & x_1^{(0)} &         0 & x_2^{(0)} & x_3^{(0)} \\
   0 & -2L_0 &          0 &           0 &    0 &         0 &         0 & x_3^{(1)} & x_4^{(1)} &         0 \\
   0 &     0 & 2x_3^{(1)} &           0 &    0 &         0 &         0 &      -L_0 &         0 & x_4^{(1)} \\
   0 &     0 &          0 &  2x_4^{(1)} &    0 &         0 &         0 &         0 &      -L_0 & x_3^{(1)} \\
 -2L_0 &   0 &          0 &           0 &    0 & x_3^{(2)} & x_4^{(2)} &         0 &         0 &         0 \\
   0 &     0 & 2x_3^{(2)} &           0 &    0 &      -L_0 &         0 &         0 &         0 & x_4^{(2)} \\
   0 &     0 &          0 &  2x_4^{(2)} &    0 &         0 &      -L_0 &         0 &         0 & x_3^{(2)}
\end{array}
\right)
.
$$
If~this matrix is of
rank~$9$
then, again, we have a degenerate~case. There~is a linear form
in~$x_1, \ldots, x_4$
only, vanishing on the three points~given. The~unique solution of the system corresponds to the square of this linear~form.

Consequently,~all the ten
\mbox{$9 \times 9$-minors}
must~vanish. Actually,~when deleting the fourth column, the corresponding minor~is
$$-16 L_0^4 L_1 L_2 (L_0x_3^{(0)} - L_1x_3^{(1)} + L_2x_3^{(2)})^2 \, . \eop$$
}
\end{lem}

\begin{rem}[{\rm Interpretation}{}]\vskip-\bigskipamount
The~relations 
established in Lemma~\ref{Groeb} may be interpreted as~follows. The~coordinates of three points of tangency form a~$3 \times 5$-matrix
$$
\left(
\begin{array}{ccccc}
0         & x_1^{(0)} & x_2^{(0)} & x_3^{(0)} & x_4^{(0)} \\
x_0^{(1)} & 0         & x_2^{(1)} & x_3^{(1)} & x_4^{(1)} \\
x_0^{(2)} & x_1^{(2)} & 0         & x_3^{(2)} & x_4^{(2)}
\end{array}
\right)
.
$$
We~may scale such that
$\smash{x_0^{(1)} = x_1^{(0)}}$
and~$\smash{x_0^{(2)} = x_2^{(0)}}$.\smallskip

\begin{iii}
\item
Then,~the leftmost
$3 \times 3$-block
is either symmetric,
i.e.,~$\smash{x_1^{(2)} = x_2^{(1)}}$,
or symmetric up to~sign.
Then,~$\smash{x_1^{(2)} = -x_2^{(1)}}$.
\item
In~the latter case, the column vector\vspace{0.4mm}
$\smash{(x_3^{(0)}, x_3^{(1)}, x_3^{(2)})^t}$
is a linear combination of the column
vectors~$\smash{(0, x_0^{(1)}, x_0^{(2)})^t}$,
$\smash{(x_1^{(0)}, 0, x_1^{(2)})^t}$,
and
$\smash{(x_2^{(0)}, x_2^{(1)}, 0)^t}$.
\end{iii}
\end{rem}

\begin{rems}
\begin{iii}
\item
In~the non-symmetric variant,\vspace{0.5mm}
$\smash{(x_4^{(0)}, x_4^{(1)}, x_4^{(2)})^t}$
is a linear combination of the column
vectors~$\smash{(0, x_0^{(1)}, x_0^{(2)})^t}$,
$\smash{(x_1^{(0)}, 0, x_1^{(2)})^t}$,
and~$\smash{(x_2^{(0)}, x_2^{(1)}, 0)^t}$,
too. The~roles of the third and fourth coordinates may be~interchanged.
\item
Actually,~in~this variant, linear dependence of the three vectors\vspace{0.3mm}
$\smash{(0, x_0^{(1)}, x_0^{(2)})^t}$,
$\smash{(x_1^{(0)}, 0, x_1^{(2)})^t}$,
and
$\smash{(x_2^{(0)}, x_2^{(1)}, 0)^t}$
is a non-trivial~condition. Observe,~they do not form a base
of~$\bbR^3$.
In~the symmetric variant, an analogous condition would be~empty.
\end{iii}
\end{rems}

\begin{rem}
For~each triple consisting of points of tangency
of~$S$
with a coordinate~hyperplane, relations of the same kind must be~fulfilled.
\end{rem}

\noindent
{\bf Proof of Theorem~\ref{main}.}
For~each of the five points of tangency, we have at least two pairs
$\{i,j\} \subset \{ 0, \ldots, 4 \}$
such
that~$x_i = x_j$.
There~are two~cases.\medskip

\noindent
{\em First case.}
Each~of the ten pairs
of~$\{ 0, \ldots, 4 \}$
appears exactly~once.\smallskip

\noindent
Without~restriction, the point of tangency
to~$H_0$
is~$(0 : 1 : 1 : t : t)$.
Again~without loss of generality,
$(1 : 0 : s : 1 : s)$
is the point of tangency
to~$H_1$.
The~structure of the remaining three points of tangency is then~fixed. The~five points form a matrix as~follows,
$$
\left(
{\arraycolsep1.5mm
\begin{array}{ccccc}
0 & 1 & 1 & t & t \\
1 & 0 & s & 1 & s \\
1 & r & 0 & r & 1 \\
t & q & t & 0 & q \\
t & t & p & p & 0 
\end{array}
}
\right)
.
$$
Lemma~\ref{Groeb} implies
that~$r = s$.
Indeed,~$r = -s$
would enforce that both
$(t,1,-s)^t$
and~$(t,s,1)^t$
are linearly dependent
of~$(0,1,1)^t$,
$(1,0,-s)^t$,
and~$(1,s,0)^t$.
This~is a contradiction
since~$(0,s-1,s+1)^t$
is not in the span of these~three.

For~the same
reason,~$p = q$.
Further,~we have
$q = \pm 1$
and~$s = \pm t$
such that we end up with four one-parameter~families,
$$
\left(
{\arraycolsep1.3mm
\begin{array}{rrrrr}
0 & 1 & 1 & t & t \\
1 & 0 & t & 1 & t \\
1 & t & 0 & t & 1 \\
t & 1 & t & 0 & 1 \\
t & t & 1 & 1 & 0 
\end{array}
}
\right)
\! , \,
\left(
{\arraycolsep0.2mm
\begin{array}{rrrrr}
0 & 1 & 1 & t & t \\
1 & 0 & -t & 1 & -t \\
1 & -t & 0 & -t & 1 \\
t & 1 & t & 0 & 1 \\
t & t & 1 & 1 & 0 
\end{array}
}
\right)
\! , \,
\left(
{\arraycolsep0.2mm
\begin{array}{rrrrr}
0 & 1 & 1 & t & t \\
1 & 0 & t & 1 & t \\
1 & t & 0 & t & 1 \\
t & -1 & t & 0 & -1 \\
t & t & -1 & -1 & 0 
\end{array}
}
\right)
\! , \,
\left(
{\arraycolsep0.2mm
\begin{array}{rrrrr}
0 & 1 & 1 & t & t \\
1 & 0 & -t & 1 & -t \\
1 & -t & 0 & -t & 1 \\
t & -1 & t & 0 & -1 \\
t & t & -1 & -1 & 0 
\end{array}
}
\right)
\! .
$$
The~linear equation
of~$S$
requires that the matrices considered are of rank at
most~$4$.
However,~in the second and third families, the determinants
$(t^2-t-1)(t^3+2t-1)$
and~$(t^2+t-1)(t^3-2t^2-1)$
have no rational~zeroes. For~the fourth family, we find
$(t+1)(t^2-t+1)(t^2+3t+1)$
for the~determinant.
But,~for~$t = -1$,
we had four equal coordinates in several of the points of~tangency. Finally,~for the first family, the determinant is
$(t+1)(t^2-3t+1)^2$
and the
value~$t = -1$
could be~possible.

The~corresponding data lead to systems of equations which are uniquely solvable up to~scaling. The~resulting quadric surface is given
by~$l = q = 0$~for
\begin{eqnarray*}
l & := & x_0 + x_1 + x_2 + x_3 + x_4 \, , \\
q & := & x_0^2 + x_1^2 - x_2^2 - x_3^2 + 3x_0x_1 + x_0x_2 - x_0x_3 - x_1x_2 + x_1x_3 - 3x_2x_3 \, .
\end{eqnarray*}
This~surface is indeed smooth and tangent to all five coordinate hyperplanes but the double
covering~$\pi\colon O \to \bP^4_\bbQ$
does not split over~it.\medskip

\noindent
{\em Second case.}
One~of the ten pairs
of~$\{ 0, \ldots, 4 \}$
appears at least~twice.\smallskip

\noindent
Without~loss of generality, the points of tangency
to~$H_0$
and~$H_1$,
respectively, are
$(0 : 1 : 1 : t : t)$
and~$(1 : 0 : 1 : s : s)$.
If~the point of tangency
to~$H_2$
were~$(1 : (-1) : 0 : 1 : (-1))$
then, by~Lemma~\ref{Groeb}, both
$(t,s,1)^t$
and~$(t,s,-1)^t$
had to be linear combinations
of~$(0,1,1)^t$,
$(1,0,-1)^t$,
and~$(1,1,0)^t$.
This~is a contradiction
since~$(0,0,2)^t$
is not in the span of these~three.

Consequently,~the five points of tangency form a matrix as~follows,
$$
\left(
{\arraycolsep0.5mm
\begin{array}{rrrrr}
0 &     1 &     1 &     t &              t \\
1 &     0 &     1 &     s &              s \\
1 &     1 &     0 &     r &              r \\
t & \pm s & \pm r &     0 &              q \\
t & \pm s & \pm r & \pm q & \phantom{\pm}0 
\end{array}
}
\right)
.
$$
Assume~that one of the
``$r$''
or~``$s$''
actually carries a minus~sign. Without~restriction, there is
``$-r$''
in the fourth~line. Then,~Lemma~\ref{Groeb} yields the contradiction that
$(t,r,q)^t$
must be a linear combination
of~$(0,1,t)^t$,
$(1,0,-r)^t$,
and~$(t,r,0)^t$.
Further,~if there were a
``$-q$''
in the fifth line then
$(1,r,r)^t$
had to be a linear combination
of~$(0,t,t)^t$,
$(t,0,-q)^t$,
and~$(t,q,0)^t$
which is not the case,~either.

Finally,~in the fourth line, we must have two pairs of equal~entries. Without~restriction, suppose that
$q = t$
and~$r = s$.
All~in all, we find a matrix of the~form
$$
\left(
{\arraycolsep1.5mm
\begin{array}{rrrrr}
0 & 1 &  1 &  t & t \\
1 & 0 &  1 &  s & s \\
1 & 1 &  0 &  s & s \\
t & s &  s &  0 & t \\
t & s &  s &  t & 0 
\end{array}
}
\right)
.
$$
For~the determinant, one
calculates~$2t^2(4s-t-1)$.
We~may conclude
that~$s = \frac{t+1}4$.

For~every~$t \neq 0$,
these data lead to systems of equations which are uniquely solvable up to~scaling. The~result is the
one-parameter~family~$S_t$
of quadric surfaces given
by~$l_t = q_t = 0$~for
\begin{eqnarray*}
l_t & := & (t-1)x_0 - 2tx_1 -2tx_2 + 2x_3 + 2x_4 \, , \\
q_t & := & (t+1)^2x_0^2 + 4(t+1)tx_1^2 + 4(t+1)tx_2^2 + 16x_3^2 - {} \\
 & & \hspace{3cm} {} - 4(t+1)tx_0x_1 - 4(t+1)tx_0x_2 + 8(t-1)x_0x_3 + {} \\
 & & \hspace{5.3cm} {} + 8(t-1)tx_1x_2 - 16tx_1x_3 - 16tx_2x_3 \, .
\end{eqnarray*}
For~each~$t \neq 0$,
the quadric
surface~$S_t$
is indeed smooth and tangent to all five coordinate~hyperplanes.

In~order to check for which values
of~$t$
the double
covering~$\pi\colon O \to \bP^4_\bbQ$
splits
over~$S_t$,
we first restrict to the
intersection~$C_t := S_t \,\cap {}``x_1 \!=\! x_0 \!+\! x_2$''.
This~is a smooth~conic for
each~$t \neq 0$.
A~parametrization
$\iota_t \colon \bP^1 \to C_t$
is given~by
{\small
\begin{eqnarray*}
&& \hspace{-0.25cm} (u : v) \mapsto
(16tu^2 :
((t^2 + 18t + 1)u^2 + 8(t + 1)tuv + 16t^2v^2) : \\
&& \hspace{-0.1cm} : ((t^2 \!+\! 2t \!+\! 1)u^2 + 8(t \!+\! 1)tuv + 16t^2v^2) :
((t^2 \!+\! 2t \!+\! 1)tu^2 + 8(t \!-\! 1)t^2uv + 16t^3v^2) : \\
&& \hspace{5.3cm} : ((t^2 + 10t + 9)tu^2 + 8(t + 3)t^2uv + 16t^3v^2)) \, .
\end{eqnarray*}}%
The~binary form
$(-3) \Delta^\prime (\iota_t (u, v))$
of degree~16 factors into
$u^6 ((t + 1)u + 4tv)^4$
and a form of degree~six which is irreducible for
general~$t$.
We~ask for the values
of~$t$
for which this sextic is a perfect~square. According~to {\tt magma}, its discriminant is equal~to
$$C (t-3) (t-1)^6 (3t - 1)^6 t^{83} (t^2 + 8t - 1)^4 (19t^3 - 82t^2 + 59t - 16)^2$$
for~$C$
a 103-digit~integer.
Over~$S_1$,
the double
covering~$\smash{\pi\colon O \to \bP^4_\bbQ}$
does not~split. The~cases
$t = 3$
and~$t = \frac13$
both yield the accumulating
subvariety~$Q$
studied in subsection~\ref{Infi}. They~are equivalent to each other under the
permutation~$(0)(13)(24)$
of~coordinates.
\eop

\end{document}